\documentclass[12pt]{article} 
\usepackage{amssymb,latexsym}
\usepackage{amsfonts}   
\input cyracc.def
 at 12truept

 at 12truept

 at 15truept

 at 12truept
 at 12truept

 at 12truept

\font\germ=eufm10 \font\eur=eusm10 
\def\qed{\hfill $\vrule height 2.5mm  width 2.5mm depth 0mm $}

\def\neweq{\setcounter{equation}{0}}
\newtheorem{theorem}{Theorem}[section]
\newtheorem{pr}[theorem]{Proposition}
\newtheorem{cor}[theorem]{Corollary}
\newtheorem{de}[theorem]{Definition}
\newtheorem{con}[theorem]{Conjecture}

\newtheorem{prb}[theorem]{Problem}
\newtheorem{ex}[theorem]{Example}

\def\g{\hbox{\germ g}}

\def\M{\hbox{\germ M}}

\def\R{\mathbb{R}}
\def\N{\mathbb{N}}
\def\Z{\mathbb{Z}}
\def\C{\mathbb{C}}

\def\Pe{\hbox{\eur P}}

\def\ds{\displaystyle}
\def\wt{\widetilde}

\def\ld{\lambda}

 \def\al{\alpha}

\def\bdoteq{\buildrel\bullet\over{=\!\!\!=}}
\begin{document}
\begin{center}
{\textbf{\Large Rigged Configurations and Catalan, \\
Stretched Parabolic Kostka Numbers 
and Polynomials : \\
Polynomiality, Unimodality and Log-concavity
}}
\end{center}

\begin{center}
\textsc{anatol n.\ kirillov}
\end{center}
\begin{center}
{\small {\it Research Institute of Mathematical Sciences ( RIMS )}} \\
{\small {\it Kyoto 606-8502, Japan }}\\ {\small}
{\it URL: ~~~~http://www.kurims.kyoto-u.ac.jp/\kern-.05cm$\tilde{\quad}$\kern-.17cm kirillov }
$$and$$
{\small {The Kavli Institute for the Physics and Mathematics of the Universe 
( IPMU ),\\ 5-1-5 Kashiwanoha,  Kashiwa, 277-8583, Japan}} \\
\vspace{8mm}

\vspace{5mm}

\vspace{8mm}

\end{center}
 
\vspace{5mm}
\begin{abstract}
We will look at the Catalan numbers from the {\it Rigged Configurations} point
 of view originated \cite{Kir} from an combinatorial  analysis of the 
Bethe Ansatz Equations associated with the higher spin anisotropic Heisenberg 
models .  Our strategy is to take a combinatorial interpretation of Catalan 
numbers $C_n$ as the number of standard Young tableaux of rectangular 
shape  $(n^2)$, or equivalently, as the Kostka number $K_{(n^2),1^{2n}}$, as 
the starting point of research. We observe  that the rectangular (or 
multidimensional)  Catalan numbers $ C(m,n)$ introduced and studied by P. 
MacMahon \cite{Mc}, \cite{Su1}, see also \cite{Su2}, can be identified with 
the Kostka number $K_{(n^m),1^{mn}}$, and therefore can be treated by Rigged 
Configurations technique. Based on this technique we study the stretched 
Kostka numbers and polynomials, and give a proof of `` a strong rationality 
`` of the stretched Kostka polynomials. This result implies a polynomiality 
property  of the stretched Kostka and stretched
  Littlewood--Richardson coefficients \cite{KT}, \cite{Ras}, \cite{Ki1}. \\

Another application of the Rigged Configuration technique presented, is 
 a new family of counterexamples to Okounkov's log-concavity conjecture 
\cite{Ok}.\\

Finally, we apply Rigged Configurations technique to give a combinatorial 
prove of the unimodality of the principal specialization of the internal 
product of Schur functions. In fact we prove a combinatorial formula for 
generalized $q$-Gaussian polynomials which is a far generalization of the 
so-called $KOH$-identity \cite{O}, as well as it manifests the unimodality 
property of the $q$-Gaussian polynomials.

\end{abstract}
%\vskip4cm

%\hrule width 3 in 
%\vspace{2mm}
{2000 Mathematics Subject Classifications: 05E05, 05E10, 05A19.}

%\vskip 0.5cm
%\pagebreak
%\begin{center}
%\textsc{Contents}
%\end{center}

%\nopagebreak
%\begin{tabular}{rrl}
%\end{tabular}

\section{Introduction}

The literature devoted to the study of Catalan
\footnote{$en.wikipedia.org/wiki/Catalan{\_}number$}
 and Narayana numbers
\footnote{~~$en.wikipedia.org/wiki/Narayana{\_}number$}
, their different combinatorial interpretations (more than $200$ in fact, \cite{St4}),  numerous generalizations, applications to Combinatorics, Algebraic 
Geometry, Probability Theory and so on and so forth, are enormous, see 
\cite{St4} and the literature quoted therein.  There exists a wide  variety 
of different  generalizations of Catalan  numbers, such as the Fuss--Catalan 
numbers
\footnote{$en.wikipedia.org/wiki/Fuss-Catalan{\_}number$}
 and  the Schr\"{o}der numbers
\footnote{wolfram.com/Schr\"{o}ederNumber.html}
, higher genus multivariable  Catalan numbers \cite{Mu},~ 
higher dimensional Catalan and Narayana numbers \cite{Mc}, \cite{Su1}, and 
many  and varied other generalizations. Each a such generalization comes from 
a generalization of a certain combinatorial interpretation of Catalan numbers,
 taken as a starting point for investigation. One a such interpretation of 
  Catalan numbers has been taken as the starting point of the present paper, 
is the well-known fact that the Catalan number $C_n$ is equal to the number of 
{\it standard} Young tableaux of the shape $(n^2)$. 

Let us look at the Catalan numbers from {\it Rigged Configurations} 
side. Since $C_n =K_{(n^2),1^{2n}}$ we can apply a {\it fermionic formula} for 
Kostka numbers \cite{Ki6}, and come
 to the following combinatorial expressions for Catalan and Narayana numbers
$$ C_n = \sum_{\nu \vdash n}~~ \prod_{j \ge 1} ~{2 n -2 ( \sum_{a \le j} \nu_{j}) + \nu_{j} -\nu_{j+1} \choose \nu_{j}-\nu_{j+1}},$$ 
where the sum runs over all partitions $\nu$ of size $n$; 
$$ N(n,k) = \sum_{\nu \vdash n \atop \nu_1 = k} ~~\prod_{j \ge 1} ~{2 n -2 ( 
\sum_{a \le j} \nu_{j}) + \nu_{j} -\nu_{j+1} \choose \nu_{j}-\nu_{j+1}},$$  
where the sum runs over all partitions $\nu$ of size $n$, ~$\nu_1=k$. \\
A $q$-versions of these formulas one can find, for example, in \cite{Ki2}. \\
Let us illustrate our fermionic formulas for $n =6$. There are $11$ partitions of size $6$. We display below the  distribution of contributions to the 
fermionic formulae for  the Catalan and Narayana numbers presented above, 
which come  from partitions $\nu$ of size $6$ and $k=1,2,\ldots,6$. \\
$N(1,6)=1$, \\
$N{2,6}={9 \choose 1} +{5 \choose 4} +1 = 9+5+1=15$, \\
$N(3,6)= {8 \choose 6} +{3 \choose 1} {7 \choose 1} +1 = 28+21+1=50$, \\
$N(4,6)= {7 \choose 4}+ {6 \choose 4}= 35+15 = 50,$\\
$N(5,6)= {6 \choose 4} =15,$ \\
$N(6,6)= 1$. \\

 A few comments in order.

$\bullet$~~$q$-versions of formulas for Catalan and Narayana numbers displayed 
above coincide with the Carlitz--Riordan $q$-analog of Catalan numbers 
\cite{St3} and $q$-analog of Narayana numbers correspondingly.

$\bullet$~~ It is well-known that partitions of $n$ with respect to the 
dominance ordering, form a lattice denoted by $L_n$. One (A.K)  can define an 
ordering
\footnote{~~This ordering is inherited from an {\it evolution} of the maximal 
configuration $\Delta(\ld,\mu)$ of type $(\ld, \mu)$ under certain 
transformations on the set of admissible configurations $C(\ld,\mu)$ coming 
from representation theory of the algebra ${\mathfrak{gl}}_n$. An evolution of the {\it maximal rigged configuration}  $(\Delta(\ld,\mu), J=\{P_{j}^{(k)}(\Delta(\ld,\mu)) \})$ induces a certain poset structure on the set of rigged 
configurations $R_{\ld,\mu}$. We expect some connections between certain 
poset (lattice ?) structures arising on the set of Rigged Configurations of 
type $(n^2,1^{2n})$ coming from  an analysis of the Bethe Ansatz, and the 
known poset and lattice structures on a variety of  sets counting by Catalan 
numbers, see e.g. \cite{St4}, \cite{Ass}. It is an interesting and important 
Problem to study poset structures on the set(s) counting by the 
 multidimensional Catalan and Narayana numbers.  Details will appear elsewhere.}
  on the set of admissible configurations of type $(\ld,\mu)$ as well. In the 
case $\ld=(n^2),~\mu = (1^{2n})$ the poset of admissible configurations of 
type $(\ld,\mu)$ is essentially the same as the lattice of partitions $L_n$. 
Therefore, to each vertex $\nu$ of the lattice $L_n$ one can attach the space 
of rigged configurations $RC_{\ld, \mu}(\nu)$  associated with partition $\nu$.
~Under a certain evolution a configuration $(\nu, J)$ evolves and touches the 
boundary of the set $RC_{\ld,\mu}(\nu)$. ~When such is the case,  ``state''~~ 
$(\nu, J)$ suffers ~``a phase transition''~, executes the wall-crossing, and 
end up as a newborn state of some space $RC_{\ld,\nu}(\nu^{'})$. A  precise 
description of this  process is the essence of the Rigged Configuration 
Bijection \cite{Ki9}, \cite{KSS}. It seems an interesting task to write out 
in full the evolution  process going on in the space of triangulations of a 
convex $(n+2)$-gon under the Rigged Configuration Bijection.

$\bullet$~~It is an open {\bf Problem} to count the number of admissible 
configurations associated with the multidimensional Catalan numbers $C(m,n)$, 
if $ m \ge 3$, and a structure of the corresponding poset, and dynamics of a 
rigged configuration in a given set of configuration. If $m=3$, the set of of 
admissible configurations consists of pairs of partitions $(\nu^{(1)},\nu^{(2)})$ such that $\nu^{(2)} \vdash n$ and $\nu^{(1)} \ge \nu^{(2)} \vee \nu^{(2)}$
\footnote{~~Recall that for any partitions $\ld$ and $\mu$, $\ld \vee \mu$ 
denotes partition corresponding to composition $(\ld_{1},\mu_{1},\ld_2,\mu_{2},
 \ldots)$.}.  One can check that the number of admissible configurations of 
type $(n^3,1^{3n})$ is equal to  $1,3,6,16,33,78$, for $n=1,2,3,4,5,6$. 

$\bullet$~~  It is well-known that the $q$-Narayana numbers
\footnote{~~Recall that the $q$- Narayana number $N(k,n ~|~q)= \frac{1-q}{1-q^n} {n \brack k}_{q}{n \brack k+1}_{q}$ .} 
 obey the symmetry 
property, namely, $N(k,n)=N(n-k+1,n)$. Therefore it implies some non trivial 
relations  among the products of $q$-binomial coefficients, combinatorial 
proofs of whose are desirable.

$\bullet$ ~~It is well-known that the Narayana number $N(k,n)$ counts the 
number of Dyck paths of the semilength $n$ with exactly $k$ peaks. Therefore, 
the set of rigged configurations $\{\nu\}$  which associated with the Catalan 
number $C_n$ and have fixed  $\nu_1=k$, is in one-to-one correspondence with the 
set of  the semilength $n$ Dyck paths with exactly $k$ peaks. 

Thus it looks natural to find and study combinatorial properties of the number of 
standard Young tableaux of an arbitrary rectangular  shape $(n^m)$,that is the 
 Kostka number  $K_{(n^m),1^{mn}}$, which  are inherent in the  classical Catalan 
and Narayana numbers. For example, one can expect that a multidimensional 
Catalan number is the sum of multidimensional Narayana ones (this is so !), or  expect that  a 
multidimensional Narayana polynomial is the $h$-vector of a certain convex 
lattice polytope, see e.g. \cite{St2} for the case of classical Catalan and 
Narayana numbers
\footnote{ The  multidimensional Catalan and Narayana numbers, as well as the 
first {\it expectatio}, had been introduced and proved by P.MacMahon \cite{Mc}. The second  {\it expectation} will be treated in the present paper, Section~3.}.

$\bullet$ Combinatorial analysis of the Bethe Ansatz Equations  \cite{Kir}, 
gives rise to a natural interpretation of the Catalan and rectangular Catalan 
and Narayana numbers in terms of rigged configurations, and pose 
{\bf Problem} to elaborate  combinatorial structures induced by rigged 
configurations on any chosen combinatorial interpretation of Catalan numbers. 
For example, how to describe  all triangulations of a convex $(n+2)$-gon 
which are in  a  ``natural''  bijection with the  set of all rigged 
configurations  ~$(\mu , J)$~ corresponding  to a given configuration  $\nu$ 
of type $((n^2),1^{2n})$ ? One can ask similar questions concerning Dyck paths and its multidimensional generalizations  \cite{Su2}, and so on. \\ 
 
In  Section~5.1 we present an example to illustrate some basic properties of 
the Rigged Configuration Bijection.\\ 
\qed

 In the present paper we are interested in to investigate  combinatorics 
related with the  higher dimensional Catalan numbers, had been introduced and 
studied in depth by P. MacMahon \cite{Mc}. It is highly possible that the 
starting point to introduce the higher dimensional Catalan numbers in 
\cite{Mc}  was  an interpretation of classical Catalan numbers as the number 
of rectangular shape $(n^2)$ standard Young tableaux mentioned above.  

Our main objective in the present paper is to look on the multidimensional 
Catalan numbers $C(m,n):=C(m,n |1)$, defined as the value of the Kostka--
Foulkes polynomials $K_{(n^m),(1^{mn})}(q)$  at $q=1$, from the point of view 
of Rigged Configurations Theory. In other words, we want to study the multidimensional 
Catalan and  Narayana numbers introduced in \cite{Mc}, \cite{Su1}, by means of a 
fermionic formula for parabolic Kostka polynomials due to the author, e.g. \cite{Ki4}, \cite{Ki1}. In particular, we apply the fermionic formula for parabolic Kostka 
polynomials cited above, to the study a stretched (parabolic) Kostka 
polynomials $K_{N \lambda,N \{\cal{R}\}}(q)$. At this way we obtain the following 
results.
\begin{theorem}~~ (Polynomiality)  ${}$ 

Let $\lambda$ be partition and $\{\cal{R}\}$ be a dominant sequence of rectangular shape partitions. ~~Then

$$ ~~~\sum_{N \ge 0} K_{N \lambda, N \cal{R}}(q)~t^N = \frac{P_{\lambda,\cal{R}}(q,t)}{Q_{\lambda, \cal{R}}(q,t)}, $$
were a polynomial $P_{\lambda,\cal{R}}(q,t)$ is such that $P_{\lambda,\cal{R}}(0,0)=1$; \\
~~a polynomial $Q_{\lambda, \cal{R}}(q,t) = \prod_{s \in S} (1- q^{s} t)$ for a some set of non-negative integers  $S:=S(\lambda,\cal{R})$, depending 
on data $\lambda$ and $\cal{R}$.
\end{theorem} 
\begin{cor}~(\cite{KT}, \cite{Ras}, \cite{Ki1}) ${}$

Let $\lambda$ be partition and $\{\cal{R}\}$ be a dominant sequence of rectangular 
shape partitions. ~~Then

$\bullet$  ~~~${\cal{K}}_{\lambda,\cal{R}}(N):=K_{N \lambda,N \cal{R}}(1)$ is 
a  polynomial of $N$ with rational coefficients.

$\bullet$~~~(Littlewood--Richardson polynomials, \cite{Mo},\cite{Ras},\cite{Ki1})

Let $\lambda$, $\mu$ and $\nu$  be partitions such that $|\lambda|+|\mu|= |\nu|$. \\
The  Littlewood--Richardson number $c_{\lambda,\mu}^{\nu}(N):= 
c_{N \lambda,N \mu}^{N \nu}$ ~is a polynomial of $N$ with rational 
coefficients.
\end{cor}
\begin{prb}  Compute
\footnote{~~ It seems that the formulas for the degree of the stretched Kostka 
polynomials stated in \cite{KT}, \cite{Ras}. \cite{Ki1} are valid only for a 
special choice of $\ld$, $\mu$ or  $\cal{R}$.} 
 the degree of polynomial ${\cal{K}}_{\lambda,,\cal{R}}(N)$.
\end{prb}
Our next objective is to define a lattice convex polytope ${\cal{P}}(n,m)$ 
which has the $h$-vector equals to the sequence of multidimensional Narayana 
numbers $\{ N(m,n;k |1),~~1 \le k \le (m-1)(n-1) \}$. 

As a preliminary step we recall the definition of a Gelfand -Tsetlin polytope. \\
Let $\lambda =(\lambda_1, \ldots, \lambda_n)$ be partition and $\mu=
(\mu_1,\ldots,\mu_n)$ be composition, $|\lambda|= |\mu|$. The Gelfand--
Tsetlin polytope  of type $(\lambda,\mu)$, denoted by $GT(\lambda, \mu)$, is 
the convex hull of all points  $(x_{ij})_{1 \le i \le j \le n} 
\in \R_{+}^{{n+1 \choose 2}}$ which satisfy the following   set of inequalities and equalities
$$  x_{i,j+1} \ge x_{ij} \ge x_{i+1,j+1} \ge 0, ~~~
x_{1j}=\lambda_{j},~~1 \le j \le n,~~~\sum_{a=1}^{j} x_{aj}= \sum_{a=1}^{j} 
\mu_{a}.$$
It is well-known that the number of integer points in the Gelfand--Tsetlin 
polytope $GT(\lambda, \mu)$, ~i.e. points $(x_{ij}) \in GT(\lambda,\mu)$ such that
~$x_{ij} \in \Z_{\ge 0},~~\forall 1 \le i \le j \le n$,  
 is equal to the Kostka number $K_{\lambda,\mu}(1)$. 
Therefore the stretched Kostka number $K_{N \lambda, N \mu}(1)$ counts the 
number of integer points in the polytope $ GT(N \lambda, N \mu) = N \cdot GT(\lambda,\mu)$. So As far as is we know, there is no general criterion to decide 
where or not the Gelfand--Tsetlin polytope $GT(\lambda,\mu)$  has only 
integral vertices, but see \cite{D},\cite{KT}, \cite{A} for particular cases treated. \\

In the present paper we are interested in the $h$-vectors of Gelfand--Tsetlin 
polytopes $GT(n,1^{d})$ and that $GT((n^k,1^{kd}), (1^k)^{n+d})$. We expect 
(cf \cite{A}) that the polytope $GT(n,1^d)$ is an integral one, but we don't 
know how  to  describe the set of  parameters $(n,k,d)$ such that  the  polytope 
$GT((n^k,(1^k)^{n+d})$ is an integral one.   
\begin{theorem}  ${}$

$(1)$~~~Let $\lambda:=\lambda_{n,d}=(n,1^{d})$~and 
$\mu =\mu_{n,d}:=(1^{n+d}).$ Then

$$\sum_{ N \ge 0} K_{N\lambda,N\mu}(q)~t^{N}= 
{C_{d,n-1}(q^{{n \choose 2}}~t,q) \over 
(q^{{n \choose 2}}~t;q))_{d(n-1)+1}},$$
where $C_{d,m}(t,q)=\sum_{k=1}^{(d-1)(m-1)}N(d,m,k~|~q)~t^{k-1}$~ stands for a 
$(q,t)$--analogue of the rectangular $(d,n)$-Catalan number. \\
In particular, the normalized volume of the Gelfand--Tsetlin polytope 
$GT((n,1^{d}),1^{n+d})$ is equal to the $d$--dimensional Catalan number
$$C_{d,n}(1,1):= (d n)!~\prod_{j=0}^{d-1}{ j! \over (n+j)! }=f^{({n^d})}=
f^{({d^n})},$$

$(2)$~~~Let $\lambda:=\lambda_{n,1,2}=(n^2,1^{2})$ ~and 
$\mu =((1,1)^{n+1}),$ $n \ge 2$. Then
$$\sum_{ N \ge 0} K_{N (n,n,1,1),N (1,1)^{n+1}}(1)~t^{N}= 
\frac{P_{2,n}(t)}{(1-t)^{4n-6}},$$
and ~$ P_{2,n}(1)=C_{n-3}~C_{n-2},$ product of two Catalan numbers. 

$(3)$ ~~~Let $\lambda:=\lambda_{n,k,d}=(n^k,1^{k d})$~and 
$\mu =((1^{k})^{n+d}),$~$d \ge 1$. ~ Then
$$\sum_{N \ge 0}~K_{N(n^k,1^{k d}).N(1^k)^{n+d}}(1)~t^N = \frac{P_{k,d,n}
(t)}{Q_{k,d,n}(t)}.$$
Moreover,  ~$P_{k,d,n}(0)=1$, ~
$$Q_{k,d,n}(t)= (1-t)^{k^2(d(n-1)-1)+2 +(k-1) \delta_{n,2} ~ \delta_{d,1})},$$
 and the polynomial 
$P_{k,d,n}(t)$ is symmetric with respect to variable $t$;
$$deg_t(P_{k,k,n}(t))= (k-1)(k(n-2)+2(\delta_{n,2}-1)).$$
\end{theorem}
One can see from Theorem~1.4, $(1)$, that the degree of the stretched Kostka 
polynomial ${\cal{K}}_{(n,1),1^{n+1}}(N):= K_{N (n,1), N( 1^{n+1})}(1)$ is 
equal to $n-1$, whereas it follows from Theorem~1.4, $(2)$ that  
$$deg_{N}({\cal{K}}_{2 (n,1), 2 (1)^{n+1}}(N)) =4 n - 7 > 3 ~ deg_{N}({\cal{K}}_{(n,1), 1^{n+1}}(N)),~~if~~N > 4.$$
Therefore one comes to an infinite family of counterexamples to Okounkov log-
concavity conjecture for the Littlewood--Richardson coefficients \cite{O}.
\begin{cor} ${}$

$\bullet$~~Let $ n \ge 3$. There exists  an integer $N_{0}(n)$ such that 
\begin{equation}
 K_{2 N (n,1), 2 N (1)^{n+1}}(1)~~ >~~ \Bigl(K_{N (n,1), N  (1)^{n+1}}(1) \Bigr)^{2}~~for~~all~~N \ge N_{0}(n);
\end{equation}
$\bullet$~~~Let $n \ge 5$ be an integer, choose ~~$\epsilon,~ 
0 \le \epsilon <  \frac{n-4}{ n - 1}$. ~~There is an integer $N_{0}(n;\epsilon)$ such that 
$$ K_{2 N (n,1), 2 N (1)^{n+1}}(1)~~ >~~ \Bigl(K_{ N (n,1), N  (1)^{n+1}}(1) \Bigr)^{3+\epsilon} ~~~~for~~ all~~~ N \ge N_{0}(n;\epsilon).$$

$(3)$ ~~~Let $n > 1+\frac{k^{2}+2}{k^{2} d}$. There exist an integer $N_{0}(n,k,d)$ such that 
$$K_{2 N (n^k,1^{kd}), 2 N (1^{k})^{n+d}}(1)~~ >~~ \Bigl(K_{ N (n^k,1^{kd}), 
N  (1^{k})^{n+d}}(1) \Bigr)^{3} ~~~~for~~ all~~~ N \ge N_{0}(n,k,d).$$
\end{cor}
For example,
$$ \bullet~~K_{2 N (5,1) ,2 N (1,1)^6}(1) > (K_{ N (5,1), N (1^6)}(1))^{3}$$
 if and only if $ N\ge 49916.$  \\
Let us recall the well-known fact that any parabolic Kostka number $K_{\lambda,
\cal{R}}(1)$ can be realized as the Littlewood--Richardson coefficient 
$c_{\lambda, M}^{\Lambda}$ for uniquely defined partitions $\Lambda$ and $M$, 
see Section~3.2 for details 
  
 ~It should be stressed that for $n=3$ the example $(1.1)$  has been  
\underline{discovered} in \cite{CDW}, and independently by 
the author (unpublished notes \cite{Ki5}). In this case the minimal value of 
$N_{0}(3)$ is equal to $23$; one can show (A.K.) that $N_{0}(4)=8$. \\

Our next objective of the present paper is to prove the unimodality of the 
principal specialization $s_{\alpha} * s_{\beta}(q,\ldots,q^{N-1})$  of Schur 
functions \cite{Ki2},\cite{Ki4}.  Proofs given in {\it loc. cit.} is based on 
an identification of the principal specialization of internal product of 
Schur functions with a certain parabolic Kostka polynomial.
\begin{theorem} ({\bf Principal specialization of the internal product of Schur functions and parabolic Kostka polynomials}) ${}$

Let $\alpha ,\beta$ be  partitions such that $|\alpha|= |\beta |$, 
 $\al_1 \le r$ and $\beta_1 \le k$. For given integer $N$ such that 
$\alpha_1+\beta_1 \le Nr$,
consider partition
$$\ld_N:=(rN-\beta_{k}^{'},rN-\beta_{k-1}^{'},\ldots ,rN-\beta_1^{'},\alpha^{'})
$$
and a sequence of rectangular shape partitions
$$R_N:=(\underbrace{(r^k),\ldots ,(r^k)}_N).$$
Then
\footnote{~~Hereinafter we shall use the notation $A(q) \bdoteq B(q)$ to mean 
that the ratio $A(q)/B(q)$ is a certain power of $q$. }

\begin{equation}
K_{\ld_N ,R_N}(q)\bdoteq s_{\al}*s_{\beta}(q,\ldots ,q^{N-1}).
\label{6.11*}
\end{equation}
\end{theorem}
Now we state a fermionic formula for polynomials 
$$V_{\alpha,\beta}^{N}(q): = s_{\alpha}*s_{\beta}(q,\ldots,q^{N-1}),$$
 which is our main tool to give a combinatorial proof of the unimodality of the 
principal specialization of the Schur functions, and that of the generalized 
$q$-Gaussian polynomials ${N \brack \lambda}_{q}$ associated with a partition 
$\lambda$, as a special case.
\begin{theorem}
Let $\alpha$ and $\beta$ be two partitions of the same size, and 
$r:=\ell(\alpha)$ be the length of $\alpha$. Then

\begin{equation} 
s_{\alpha}*s_{\beta}(q,\ldots,q^{N-1})= \sum_{\{\nu \}} q^{c(\{\nu\})} \prod_{k,j \ge 1} { P_{j}^{(k)}(\nu)+m_j(\nu^{(k)})+N(k-1)\delta_{j,\beta_{1}} \theta(r-k) \brack P_j^{(k)}(\nu)}_{q},
\end{equation}
where the sum runs over the set of admissible configurations $\{\nu \}$ of 
type $([\alpha,\beta]_{N}, (\beta_1)^N)$. Here for any partition $\lambda$, 
~$\lambda_j$ denotes its $j$-th component.
\end{theorem}
See Section~4, Theorem~4.2 for details concerning notation.
An important property which is specific to admissible configurations of type 
$[\alpha,\beta]_N,\beta_1^n)$, is the following relations
$$2c(\nu)+\sum_{k,j \ge 1} P_{j}^{(k)}(\nu)\biggl[m_j(\nu^{(k)})+N(k-1) 
\delta_{j,\beta_{1}} \theta(r-k) \biggr] = N |\alpha|,$$
which imply the unimodality of polynomials $V_{\alpha,\beta}^{N}(q)$, and 
${N \brack \alpha}_q = V_{\alpha, (|\alpha|)}^{N}(q)$.  Let us \underline{stress} 
that the sum in the $RHS(1.3)$ runs over the set of admissible configurations of 
type $([\alpha,\beta]_N ,(\beta_1)^{N})$. Remark, that the $RHS(1.3)$ has a 
natural generalization to the case $|\al| \equiv |\beta|~(mod~N)$, but in this case a
representation-theoretical meaning of the $LHS(1.3)$ is unclear to the author.\\

{\bf Acknowledgments} ${}$

The present paper is an extended and update version of notes have been 
paved the way to my Colloquium talk at the University of Queensland, Brisbane, 
Australia, November 2013, and is based 
on unpublished preprint \cite{Ki5}. I'm very appreciate to Professor Ole Warnaar (University of Queensland) for kind invitation, 
fruitful discussions and financial support of my visit.

\section{Higher dimensional Catalan and Narayana numbers, 
\cite{Mc}, \cite{Ki1}, \cite{Su1}}

\subsection { Rectangular Catalan and Narayana polynomials,
and MacMahon polytope, \cite{Ki1} }

\subsubsection{Rectangular Catalan and Narayana numbers and polynomials} 
Define {\bf rectangular Catalan polynomial}
\begin{equation}
C(n,m|q)=\frac{(q;q)_{nm}}{\ds\prod_{i=1}^n\prod_{j=1}^m(1-q^{i+j-1})}=
 [d~n]_q !~\prod_{j=0}^{d-1}~\frac{[j]_q !}{[n+j]_q !},
\label{2.11}
\end{equation}
where  $[n]_q:=\frac{1-q^n}{1-q}$ stands for the $q$-analogue of an integer $n$, and by definition $[n]_q !:= \prod_{j=1}^{n}~[j]_q$. 
\begin{pr}
\begin{equation}
q^{m\scriptsize{\left(\begin{array}{c}n\\ 2\end{array}\right)}}C(n,m|q)=
K_{(n^m),(1^{nm})}(q). \label{2.12}
\end{equation}
\end{pr}
Thus, $C(n,m|q)$ is a polynomial of degree $nm(n-1)(m-1)/2$ in the
variable $q$ with non--negative integer coefficients. Moreover,
$$C(n,2|q)=C(2,n|q)=c_n(q)=\frac{1-q}{1-q^{n+1}}\left[\begin{array}{c}
2n\\ n\end{array}\right]_q
$$
coincides with "the most obvious" $q$--analog of the Catalan numbers,
%$C_n$ defined by Carlitz~\cite{Car}.
see e.g. \cite{Fur}, p.255, or \cite{St3}, ~and \cite{Mc},
$$C(n,3 |q) =\frac{[2]_q~[3~n]_q~ !}{[n]_q~!~[n+1]_q~!~[n+2]_q~!}.$$

 It follows from (\ref{2.12}) that the rectangular Catalan
number $C(n,m|1)$ counts the number of {\it lattice} words
$$w=a_1a_2\cdots a_{nm}$$
of weight $(m^n),$ i.e. lattice words in which each $i$ between 1
and $m$ occurs exactly $n$ times. Let us recall that a word
$a_1\cdots a_p$ in the symbols $1,\ldots ,m$ is said to be a {\it
lattice} word, if for $1\le r\le p$ and $1\le j\le m-1$, the
number of occurrences of the symbol $j$ in $a_1\cdots a_r$ is not
less than the number of occurrences of $j+1$:
\begin{equation}
\#\{ i|1\le i\le r~~{\rm and}~~a_i=j\}\ge\#\{ i|1\le i\le r~~{\rm
and}~~a_i=j+1\}. \label{2.12a}
\end{equation}

For any word $w=a_1\cdots a_k$, in which each $a_i$ is a positive
integer, define the major index
$${\rm maj}(w)=\sum_{i=1}^{k-1}i\chi (a_i>a_{i+1}),
$$
and the number of descents~~~
$${\rm des}(w)=\ds\sum_{i=1}^{k-1}\chi (a_i>a_{i+1}).$$

Finally, for any integer $k$ between 0 and $(n-1)(m-1)$, define {\bf
rectangular $q$--Narayana} number
$$N(n,m;k~|~q)=\sum_wq^{{\rm maj}(w)}, $$
where $w$ ranges over all lattice words of weight $(m^n)$ such that ${\rm
des}(w)=k$. 

 Equivalently, $N(n,m;k)$ is equal to the number of rectangular standard 
Young tableaux with $n$ rows and $m$ columns having $k$ descents, i.e.
$k$ occurrences of an integer $j$ appearing in a lower row that that $j+1.$

\begin{ex} {\rm Take $n=4$, $m=3$, then
$$\sum_{k=0}^6 N(3,4;k~|~1)t^k=1+22t+113t^2+190t^3+113t^4+22t^5+t^6.
$$

%This example shows that if $n,m \ge 3$, it is unlikely that there exists a
%simple combinatorial formula for the rectangular Narayana numbers
%$N(n,m;k~|~1)$, but see formula (\ref{2.12b}).
}
\end{ex}
We summarize the basic known properties of the rectangular Catalan and 
Narayana numbers in Proposition~2.3  below.
\begin{pr} ${}$  (\cite{Mc}, \cite{Su1}, \cite{Ki4})

$({\bf A})$ ({\bf Lattice words and rectangular Catalan numbers}) ${}$

$C(n,m|q)=\ds\sum_wq^{{\rm maj}(w)}$, ~~where $w$
ranges over all lattice words of weight $(m^n)$;
\vskip 0.3cm
%$({\bf B})$ ({\bf recurrence relation for multiple Narayana numbers, $q=1$ 
%cas%e})

\vskip 0.3cm
$({\bf B})$ ({\bf Bosonic formula for multidimensional Narayana numbers})
\begin{equation}
N(n,m;k~|~q)=\sum_{a=0}^{k}(-1)^{k-a}~q^{{k-a \choose 2}}~
\left[\begin{array}{c}n~m+1\\ k-a\end{array}\right]_q~
\prod_{b=0}^{n-1}{[b]!~[m+a+b]! \over [m+b]!~[a+b]!}, \label{2.13b}
\end{equation}
\vskip 0.3cm
$({\bf C})$ ({\bf Summation formula}) Let $r$ be a positive integer, then
$$ \sum_{k=0}^{r} \left[\begin{array}{c}n~m+r-k\\ r-k\end{array}\right]_q~
N(n,m;k~|~q)=\prod_{a=0}^{m-1}{[a]!~[n+r+a]! \over [n+a]!~[r+a]!}=$$
$$= \prod_{a=0}^{n-1}{[a]!~[m+r+a]! \over [m+a]!~[r+a]!}=
\prod_{a=0}^{r-1}{[a]!~[n+m+a]! \over [n+a]!~[m+a]!}.$$
\vskip 0.3cm
$({\bf D})$ ({\bf Symmetry})
$$N(n,m;k~|~q)=q^{nm((n-1)(m-1)/2-k)}N(n,m;(n-1)(m-1)-k~|~q^{-1})=N(m,n;k~|~q),$$ 
for any integer $k$, ~$0 \le k \le (n-1)(m-1)/2$;
\vskip 0.3cm
$({\bf E})$ ~~ ({\bf $q$-Narayana numbers}) 
$$N(2,n;k~\vert~q)=q^{k(k+1)}\ds\frac{1-q}{1-q^n}~{n \brack k}_{q}~{n \brack k+1}_{q}  \bdoteq {\rm dim}_qV_{(k,k)}^{\g
l(n-k+1)},~~ 0\le k\le n-1, $$ 
where $V_{(k,k)}^{{\g}l(n-k+1)}$ stands for the irreducible
representation of the Lie algebra ${\g} l(n-k+1)$ corresponding
to the two row partition $(k,k)$; recall that for any finite
dimensional $\g l(N)$--module $V$ the symbol ${\rm dim}_qV$
denotes its $q$--dimension, i.e. the principal specialization of
the character of the module $V$: \vskip -0.5cm
$${\rm dim}_qV=({\rm ch}V)(1,q,\ldots ,q^{N-1});
$$
\vskip 0.3cm
$({\bf F})$ $N(n,m;1~|~1)=\ds\sum_{j\ge 2}\left(\begin{array}{c}n\\
j\end{array}\right)\left(\begin{array}{c}m\\ j\end{array}\right)
=\left(\begin{array}{c}n+m\\ n\end{array}\right)-nm-1$;
\vskip 0.3cm
$({\bf G})$ ({\bf Fermionic formula for $q$--Narayana numbers}, \cite{Ki1})

\begin{equation}
q^{m\scriptsize\left(\begin{array}{c}n\\ 2\end{array}\right)}
N(n,m;l~|~q)=\sum_{\{\nu\}}q^{c(\nu)}\prod_{k,j\ge 1}\left[
\begin{array}{c}P_j^{(k)}(\nu)+m_j(\nu^{(k)})\\ m_j(\nu^{(k)})\end{array}
\right]_q, \label{2.12b}
\end{equation}
summed over all sequences of partitions $\{\nu\} =\{\nu^{(1)},\nu^{(2)},
\ldots ,\nu^{(m-1)}\}$ such that

$\bullet$ $|\nu^{(k)}|=(m-k)n$, $1\le k\le m-1$;

$\bullet$ $(\nu^{(1)})_1'=(m-1)n-l$, i.e. the length of the first column
of the diagram $\nu^{(1)}$ is equal to $(m-1)n-l$,~$l= 0, \ldots,(m-1)(n-1)$;

$\bullet$ $P_j^{(k)}(\nu):=Q_j(\nu^{(k-1)})-2Q_j(\nu^{(k)})+
Q_j(\nu^{(k+1)})\ge 0$, for all $k,j\ge 1$,\\
where by definition we put
$\nu^{(0)}=(1^{nm})$; for any diagram $\ld$ the number
$Q_j(\ld)=\ld_1'+\cdots\ld_j'$ is equal to the number of cells in the
first $j$ columns of the diagram $\ld$, and $m_j(\ld)$ is equal to the
number of parts of $\ld$ of size $j$;

$\bullet$ $c(\nu)=\ds\sum_{k,j\ge
1}\left(\begin{array}{c}(\nu^{(k-1)})_j' -(\nu^{(k)})_j'\\
2\end{array}\right)$.
\end{pr}
\begin{ex} Consider the case $m=3,~n=4$. In this case $C(3,4~|~1)=462$, and the 
sequences of Narayana numbers is $(1,22,113,190,113,22,1)$. Let us display 
below the distribution of  Narayana numbers which is coming from  the counting
 the number of admissible rigged configurations of type $((4^3),(1^{12}))$ 
according to the number $(m-1) n-\ell(\nu^{(1)})$, where $\ell(\nu^{(1)})$  
denotes the length of the first configuration $\nu^{(1)}$: \\
$N(3,4;0~ |~1)=1$,~~ $N(3,4;1~ |~1)=1 +21$,~~$N(3,4;2~ |~1)=15+35+63$,\\
$N(3,4;3~ |~1)=140+15+35$, ~~$N(3,4;4~ |~1)=21+28+63$,\\~~$N(3,4;5~ |~1)=6+16$,~~
$N(3,4;6~ |~1)=1$.
\end{ex}
\begin{con} If $1\le k\le (n-1)(m-1)/2$, then
$$N(n,m;k-1~|~1)\le N(n,m;k~|~1),
$$
i.e. the sequence of rectangular Narayana numbers
$\{N(n,m;k~|~1)\}_{k=0}^{(n-1)(m-1)}$ is symmetric and {\bf unimodal}.
\end{con}

For definition of unimodal polynomials/sequences see e.g.
\cite{St2}, where one may find a big variety of examples of
unimodal sequences which frequently appear in Algebra,
Combinatorics and Geometry.

\subsubsection{ Volume of the MacMahon polytope and rectangular Catalan and Narayana numbers }
 ~~Let ${\M}_{mn}$ be the 
convex polytope in $\R^{nm}$ of all points ${\mathbf x}=(x_{ij})_{1\le i\le n,
1\le j\le m}$ satisfying the following conditions
\begin{equation}
0\le x_{ij}\le 1, ~~~x_{ij}\ge x_{i-1,j},~~~x_{ij}\ge x_{i,j-1},
\label{2.14a}
\end{equation}
for all pairs of integers $(i,j)$ such that $1\le i\le n$, $1\le
j\le m$, and where by definition we set $x_{i0}=0=x_{0j}$.

We will call the polytope $\M_{nm}$ by {\it MacMahon polytope}.
The MacMahon polytope is an integral polytope of dimension $nm$
with $\left(\begin{array}{c}m+n\\ n\end{array}\right)$ vertices
which correspond to the set of (0,1)--matrices satisfying
(\ref{2.14a}).

If $k$ is a positive integer, define $i(\M_{nm};k)$ to be the number of
points ${\bf x}\in\M_{nm}$ such that $k{\mathbf x}\in\Z^{nm}$. Thus,
$i(\M_{nm};k)$ is equal to the number of plane partitions of rectangular
shape $(n^m)$ with all parts do not exceed $k$.  By a theorem of
MacMahon (see e.g. \cite{Ma}, Chapter~I, \S 5, Example~13)
\begin{equation}
i(\M_{nm};k)=\prod_{i=1}^n\prod_{j=1}^m\frac{k+i+j-1}{i+j-1}.
\label{2.14b}
\end{equation}
It follows from (\ref{2.14b}) that the Ehrhart polynomial ${\cal
E}(\M_{nm};t)$ of the MacMahon polytope $\M_{nm}$ is completely resolved
into linear factors:
$${\cal E}(\M_{nm};t)=\prod_{i=1}^n\prod_{j=1}^m\frac{t+i+j-1}{i+j-1}.
$$
Hence, the normalized volume
$${\wt{\rm vol}}(\M_{nm})=(nm)!{\rm vol}(\M_{nm})$$
of the MacMahon polytope $\M_{nm}$ is equal to the rectangular
Catalan number $C(n,m|1)$, i.e. the number of standard Young
tableaux of  rectangular shape $(n^m)$. We refer the reader to
\cite{St3}, Section~4.6, and \cite{Hi}, Chapter~IX, for
definition and basic properties of the Ehrhart polynomial ${\cal
E}({\Pe};t)$ of a convex integral polytope ${\Pe}$.

\begin{pr}
\begin{equation}
\sum_{k\ge0}i(\M_{nm};k)z^k=\left(\sum_{j=0}^{(n-1)(m-1)}N(n,m;j)
z^j\right)/(1-z)^{nm+1}, \label{2.14c}
\end{equation}
where
$$N(n,m;j):=N(n,m;j|1)$$
denotes the rectangular Narayana number. 
\end{pr}
Thus, the sequence of Narayana numbers
$$(1=N(n,m;0),N(n,m;1),\ldots ,N(n,m;(n-1)(m-1))=1)
$$
is the {\it $h$--vector} (see e.g. \cite{St3}, p.~235) of the
MacMahon polytope. In the case $n=2$ (or $m=2$) all these results
may be found in \cite{St3}, Chapter~6, Exercise~6.31.\\ 
{\bf Question}.~~({\it Higher associahedron}) Does there exist an
$(m-1)(n-1)$--dimensional integral convex (simplicial?) polytope
$Q_{n,m}$ which has $h$--vector
$$h=(h_0(Q_{n,m}),h_1(Q_{n,m}),\ldots ,h_{(n-1)(m-1)}(Q_{n,m}))$$
given by the {\it rectangular Narayana numbers} $N(n,m;k)$ :
$$\sum_{i=0}^{(n-1)(m-1)}h_i(Q_{n,m})t^{i}=C(n,m|t)~?
$$
We refer the reader to \cite{Hi}, Chapter~I, \S 6 and Chapter~III,
for definitions and basic properties of the $h$--vector of a
simplicial polytope; see also, R.~Stanley (J. Pure and Appl.
Algebra {\bf 71} (1991), 319-331).

An answer on this question is known if either $n$ or $m$ is
equal to 2, see e.g. R.~Simion (Adv. in Appl. Math. {\bf 18}
(1997), 149-180, Example~{\bf 4} (the Associahedron)). 
%Note finally, that MacMahon's polytope and Narayana's  numbers will
%appear again in Section~\ref{kfs}, Exercise~3{\bf e}.

\begin{de}(\cite{Mc}, \cite{Su1})~~Define {\bf rectangular} Schr\"oder 
polynomial
$$S(n,m|t):=C(n,m|1+t),$$
and put
$$S(n,m|t)=\ds\sum_{k\ge0}^{(n-1)(m-1)}S(n,m||k)t^k.
$$
\end{de}
 A  combinatorial interpretation  of the numbers
$S(n,m||k)$ and $S(n,m|1)$ has been done by R.~Sulanke \cite{Su1}.

\subsubsection { Rectangular Narayana and Catalan numbers, 
and $d$ dimensional lattice paths, \cite{Su1}}

Let ${\cal C}(d,n)$ denote the set of $d$--dimensional lattice paths using 
the steps
$$X_1=(1,0,\cdots,0),~X_2=(0,1,\cdots,0),\cdots,X_d=(0,0,\cdots,1),$$
running from $(0,0,\cdots,0)$ to $(n,n,\cdots,n)$, and lying in the region
$$ \{ (x_1,x_2,\cdots,x_d) \in \R_{\ge 0}^{d} ~\vert~ 
x_1 \le x_2 \le \cdots \le x_d  \}.$$
For each path $P:= p_1p_2\cdots p_{nd} \in {\cal C}(d,n)$ define the statistics
$$asc(P):= \# \{j ~ \vert~p_jp_{j+1}= X_k~X_l ,~k < l \}.$$
\begin{de} The $n$--th ~$d$--dimensional MacMahon--Narayana number of level 
$k$, ~$MN(d,n,k)$ counts the paths $P \in {\cal C}(d,n)$  with $asc(P)=k.$
\end{de}
\begin{pr}( Cf \cite{Su1} )
For any $d \ge 2$ ~and for ~$0 \le k \le (d-1)(n-1),$
$$MN(d,n,k)=\sum_{j=0}^{k}(-1)^{k-j}{dn+1 \choose k-j}~
\prod_{a=0}^{j-1}{a!~(d+n+a)! \over (d+a)!~(n+a)! }. $$
\end{pr} 
Note that the product $\prod_{a=0}^{j-1}{a!~(d+n+a)! \over (d+a)!~(n+a)! }$ is
equal to the number of plane partitions of the rectangular shape $(n^d)$,
all the parts do not exceed $j.$
\begin{de}
For $d \ge 3$ and $n \ge 1$ the $n$-th $d$-Narayana polynomial defined to be
$$N_{d,n}(t)= \sum_{k =0}^{(d-1)(n-1)} MN(d,n,k)~t^k.$$
\end{de}
\begin{cor} (Recurrence relations, \cite{Su1}) For any integer $m \ge 0$ one 
has
$$\sum_{k=0}^{m}{dn+m-k \choose m-k}~MN(d,n,k)=
\prod_{a=0}^{d-1}{a!~(n+m+a)! \over (n+a)!~(m+a)! }.$$
\end{cor}
\begin{cor} The MacMahon--Narayana number $MN(d,n,k)$ is equal to the
rectangular Narayana number $N(d,n;k).$ 
\end{cor}

\subsubsection{Gelfand--Tsetlin polytope ~~$GT((n,1^{d}),(1)^{n+d})$~~ and 
rectangular Narayana numbers}

\begin{theorem} ~~~Let $\lambda:=\lambda_{n,d}=(n,1^{d})$~and 
$\mu =\mu_{n,d}:=(1^{n+d}).$ Then
$$\sum_{ N \ge 0} K_{N\lambda,N\mu}(q)~t^{N}= 
{C_{d,n-1}(q^{{n \choose 2}}~t,q) \over 
(q^{{n \choose 2}}~t;q))_{d(n-1)+1}},$$
where $C_{d,m}(t,q)=\sum_{k=0}^{(d-1)(m-1)}N(d,m,k~|~q)~t^k$ stands for a 
$(q,t)$--analog of the rectangular $(d,n)$-Catalan number. \\
In particular, the normalized 
volume of the Gelfand--Tsetlin polytope 
$GT((n,1^{d}),1^{n+d})$ is equal to the $d$--dimensional Catalan number
$$C_{d,n}(1,1):= (dn)!~\prod_{j=0}^{d-1}{ j! \over (n+j)! }=f^{({n^d})}=
f^{({d^n})},$$
where for any partition $\lambda$, $f^{\lambda}$ denotes the number of standard
Young tableaux of shape $\lambda.$ \\
\end{theorem}

\section{Rigged configurations, stretched Kostka numbers, log-concavity and 
unimodality}
\subsection{Stretched Kostka numbers $K_{N(n^k.1^{kd}),N(1^k)^{n+d}}(1)$}
\begin{theorem}
$(1)$
$$\sum_{ N \ge 0} K_{N (n,n,1,1)),N((1,1)^{n+1})}(1)~t^{N}= 
\frac{P_{2,n}(t)}{(1-t)^{4n-6}},$$
and ~$ P_{2,n}(1)=C_{n-3}~C_{n-2}.$ 

$(2)$ Let $d \ge 1$, then
$$\sum_{N \ge 0}~K_{N(n^k,1^{k d}).N(1^k)^{n+d}}(1)~t^N = \frac{P_{k,d,n}
(t)}{Q_{k,d,n}(t)}.$$
Moreover,  ~$P_{k,d,n}(0)=1$, ~
$$Q_{k,d,n}(t)= (1-t)^{k^2(d(n-1)-1)+2 +(k-1) \delta_{n,2} ~ \delta_{d,1})},$$
 and the polynomial 
$P_{k,d,n}(t)$ is symmetric with respect to variable $t$;
$$deg_t(P_{k,k,n}(t))= (k-1)(k(n-2)+2(\delta_{n,2}-1)).$$
\end{theorem}
\qed 

For example, assume that  $d=1$ and set $P_{k,n}(t):=P_{k,1,n}(t).$~~Then \\
$P_{2,3}(t)=1,~P_{2,4}(t)=(1,0,1),
P_{2,5}(t)=(1,1,6,1,1),~~ P_{2,6}(t)=(1,3,21,20,21,3,1),$ \\
$P_{2,7}(t)=(1,6,56,126,210,126,56,6,1),$ \\
$P_{2,8}(t)=(1,10,125,500,1310,1652,1310,500,125,10,1)$,
~$P_{3,3}(t)=(1,-1,1)$, \\~
$P_{3,4}(t)=(1,0,20,20,55,20,20,0,1),P_{3,5}(t)=(1,6,141,931,4816,13916,$ \\
$27531,33391,27531,
13916,4816,931,141,6,1), P_{4,1,3}(t)=(1,-3,9,-8,9,-3,1) =P_{4,2,2}(t)$.\\

It follows from the duality theorem  for parabolic Kostka 
polynomials \cite{Ki1}  that 
$$K_{N\lambda,N\mu}(1)=K_{(Nn,N^{d})^{'},((N)^{n+d})^{'})}(1)=
K_{((d+1)^{N},1^{N(n-1)}),(1^{N})^{n+d})}(1),$$
and 
$$K_{((2d+2)^{N},2^{N(n-1)}),((2)^{N})^{n+d})}(1)=
K_{(Nn,Nn,N^{2d}),((N,N)^{n+d})}(1).$$
Now consider the case $d=1$, that is  $\lambda=(n,1),~\mu=(1^{n+1})$.~ Then
$$K_{N \lambda,N \mu}(1) = K_{(Nn,N),(N^{n+1})}(1) = {N+n-1 \choose n-1}.$$
The second equality follows from a more general result \cite{Ki2}, \cite{Ki1},
\begin{pr} Let $\lambda$ be a partition and $N$ be a positive integer. ~
Consider partitions $\lambda_{N} := (N |\lambda|, \lambda)$ and $\mu_N:= 
(|\lambda|^{N+1})=(\underbrace{|\lambda|,\ldots,|\lambda||}_{N+1})$. Then
$$K_{\lambda_N,\mu_N}(q)\bdoteq {N \brack \lambda} = dim_q V_{\lambda}^{{\mathfrak{gl}}(N)},$$
where the symbol ~$P(q) \bdoteq R(q)$ means that the ratio $P(q)/R(q)$ is a power of $q$; the symbol ${N \brack \lambda}$ stands for the generalized Gaussian coefficient corresponding to a partition $\lambda$, see \cite{Ma} for example.  \end{pr}
\subsection{Counterexamples to Okounkov's log-concavity conjecture}

On the other hand,
$$K_{2 \lambda_{N}, 2 \mu_{N}}(1) = K_{N(n,n,1,1),N(1,1)^{n+1}}(1)= 
Coeff_{t^{N}} \Biggl(\frac{P_{2,n}(t)}{(1-t)^{4n-6}} \Biggr) .$$
Therefore the number $K_{2 \lambda_{N}, 2 \mu_{N}}(1)$ is a {\bf polynomial} 
of the degree $4~n-7$ with respect to parameter $N$. Recall that the number 
$K_{N \lambda,N \mu}(1)={N+n-1 \choose n-1}$ is a polynomial of degree $n-1$ 
with respect to parameter $N$. Therefore we come to the following infinite 
set of examples which violate the log-concavity Conjecture stated by A. 
Okounkov \cite{Ok}. 
\begin{cor} For any integer $n > 4$, there exists a constant $N_{0}(n)$ such 
that 
$$K_{2 \lambda_{N}, 2 \mu_{N}}(1) > \Bigl(K_{N \lambda,N \mu}(1) \Bigr)^{3},$$
for all $N > N_{0}(n)$.
\end{cor}
Recall that $\lambda=(n,1),~\mu=(1^{n+1})$. \\

\underline{Now take $n=3$.}~One has \cite{CDW}
$$K_{N(3,1),N(1^4)}(1) ={N+2 \choose 2},~~K_{N(3,3,1,1),N(1,1)^4}(1)=
{N+5 \choose 5}.$$
One can check \cite{CDW} that 
$$K_{N(3,3,1,1),N(1,1)^4}(1) > \Bigl(K_{N(3,1),N(1^3)}(1) \Bigr)^2. $$
if (and only if) $N \ge 21.$

Indeed,
$$K_{N(3,3,1,1),N(1,1)^4}(1)-\Bigl(K_{N(3,1),N(1^3)}(1) \Bigr)^2=
\frac{N^2-18 N-43}{20}~{n+2 \choose 3}..$$

\underline{Now take $n=4$.} ~One has
$$K_{N(4,1),N(1^5)}(1) ={N+3 \choose 3},~~K_{N(4,4,1,1),N(1,1)^5}(1)=
{N+9 \choose 9} + {N+7 \choose 9}.$$
One can check that 
$$K_{N(4,4,1,1),N(1,1)^5}(1) > \Bigl(K_{N(4,1),N(1^5)}(1) \Bigr)^2, $$
if (and only if ) $N \ge 8$, \\
\underline{Now take $n=5$.}

\begin{pr} Let $\nu_N:=N(5,1)$~and $\eta_N:= N(1)^{6}.$~~Then 
$$\bullet~~K_{2\nu_N,2\eta_N}(1) > (K_{\nu_N,\eta_N}(1))^{2}$$
 if and only if $ N\ge 6,$   
$$ \bullet~~K_{2\nu_N,2\eta_N}(1) > (K_{\nu_N,\eta_N}(1))^{3}$$
 if and only if $ N\ge 49916.$   
\end{pr}
Indeed,
$$K_{N(5,1),N(1^6)}(1) ={N+4 \choose 4},$$
$$~~K_{N(5,5,1,1),N(1,1)^6}(1)= {N+13 \choose 13} + {N+12 \choose 13}+ 
6 {N+11 \choose 13 } + {N+10 \choose 13}+ {N+9 \choose 13},$$
and ~~$ 51891840 \times~\Bigl[ K_{N(5,5,1,1),N(1,1)^6}(1) - \Bigl(K_{N(5,1),N(1^6)}(1)\Bigr)^3 \Bigr]= {N+4 \choose 5}~ \times$ \\
$(-78631416 - 172503780~ N - 174033932~ N^2 - 101206400~ N^3 - 35852065~ N^4 -$
 
 $ 7638110~ N^5 - 899548~ N^6 - 44990~ N^7 + N^8)$. \\

Note, see e.g. \cite{Ma}, that for any set of partitions $\lambda, \mu^{(1)}, \ldots,\mu^{(p)}$ the parabolic Kostka number $K_{\lambda,\mu^{(1)},\ldots,\mu^{(p)}}(1)$ is equal to the Littlewood--Richardson number $c_{\lambda,M}^{\Lambda}$, where partitions $\Lambda \supset M$ are such that $\Lambda \setminus M = 
 \coprod_{i} \mu^{(i)}$ is a disjoint union of partitions $\mu^{(i)},~i=1,\ldots,p$.

\section{Internal product of Schur functions}
\label{ipsf}
\neweq

The irreducible characters $\chi^{\ld}$ of the symmetric group $S_n$ are
indexed in a natural way by partitions $\ld$ of $n$. If $w\in S_n$, then
define $\rho (w)$ to be the partition of $n$ whose parts are the cycle
lengths of $w$. For any partition $\ld$ of $m$ of length $l$, define the
power--sum symmetric function
$$p_{\ld}=p_{\ld_1}\ldots p_{\ld_l},$$
where $p_n(x)=\sum x_i^n$. For brevity write $p_w:=p_{\rho (w)}$.
The Schur functions $s_{\ld}$ and power--sums $p_{\mu}$ are
related by a famous result of Frobenius
\begin{equation}
s_{\ld}=\frac{1}{n!}\sum_{w\in S_n}\chi^{\ld}(w)p_w. \label{6.1}
\end{equation}
For a pair of partitions $\al$ and $\beta$, $|\al|=|\beta|=n$, let us
define the internal product  $s_{\al}*s_{\beta}$ of Schur functions
$s_{\al}$ and $s_{\beta}$:
\begin{equation}
s_{\al}*s_{\beta}=\frac{1}{n!}\sum_{w\in S_n}\chi^{\al}(w)\chi^{\beta}(w)p_w.
\label{6.2}
\end{equation}
It is well--known that
$$s_{\al}*s_{(n)}=s_{\al},~~ s_{\al}*s_{(1^n)}=s_{\al'},
$$
where $\al'$ denotes the conjugate partition to $\al$.

Let $\al ,\beta ,\gamma$ be partitions of a natural number $n\ge 1$,
consider the following numbers
\begin{equation}
g_{\al\beta\gamma}=\frac{1}{n!}\sum_{w\in S_n} \chi^{\al}(w)
\chi^{\beta}(w)\chi^{\gamma}(w). \label{6.3}
\end{equation}
The numbers $g_{\al\beta\gamma}$ coincide with the structural constants
for multiplication of the characters $\chi^{\al}$ of the symmetric group
$S_n$:
\begin{equation}
\chi^{\al}\chi^{\beta}=\sum_{\gamma}g_{\al\beta\gamma}\chi^{\gamma}.
\label{6.4}
\end{equation}
Hence, $g_{\al\beta\gamma}$ are non--negative integers. It is clear that
\begin{equation}
s_{\al}*s_{\beta}=\sum_{\gamma}g_{\al\beta\gamma}s_{\gamma}. \label{6.5}
\end{equation}

\subsection{Internal product of Schur functions, principal specialization,
fermionic formulas and  unimodality} 

Let $N\ge 2$, consider the principal specialization $x_i=q^i$, $1\le i\le
N-1$, and $x_i=0$, if $i\ge N$, of the internal product of Schur functions
$s_{\al}$ and $s_{\beta}$:
\begin{equation}
s_{\al}*s_{\beta}(q,q^2,\ldots ,q^{N-1})=\frac{1}{n!}\sum_{w\in S_n}
\chi^{\al}(w)\chi^{\beta}(w)\prod_{k\ge 1}\left(
\frac{q^k-q^{kN}}{1-q^k}\right)^{\rho_k(w)}, \label{6.9a}
\end{equation}
where $\rho_{k}(w)$ denotes the number of the length $k$ cycles of $w$.\\

By a result of R.-K.~Brylinski \cite{Br2}, Corollary~5.3, the
polynomials
$$s_{\al}*s_{\beta}(q,\ldots ,q^{N-1})
$$
admit the following interpretation. Let $P_{n,N}$ denote the
variety of $n$ by $n$ complex matrices $z$ such that $z^N=0$.
Denote by
$$R_{n,N}:=\C [P_{n,N}]$$
the coordinate ring of polynomial functions on $P_{n,N}$ with
values in the field of complex numbers $\C$. This is a graded
ring:
$$R_{n,N}=\ds\oplus_{k\ge 0}R_{n,N}^{(k)},$$
where $R_{n,N}^{(k)}$ is a finite dimensional $\g l(n)$--module
with respect to the adjoint action. Let $\al$ and $\beta$ be
partitions of common size. Then \cite{Br2}
$$s_{\al}*s_{\beta}(q,\ldots ,q^{N-1})=\sum_{k\ge 0}\langle V_{[\al ,
\beta]_n},P_{n,N}^{(k)}\rangle q^k,
$$
as long as $n\ge\max (Nl(\al),Nl(\beta),l(\al)+l(\beta))$. Here the symbol 
$\langle \bullet~,\bullet ~\rangle$ denotes the scalar product on the ring of 
symmetric functions such that $\langle s_{\lambda}~,s_{\mu} ~\rangle = \delta_{\lambda,\mu}.$ ~In other words, if $V$ and $W$ be two $GL(N)$-modules, then  
$\langle V , W  \rangle = \dim Hom_{GL(N)}(V,W)$. \\

 Let us remind below one of the main result obtained in  \cite{Ki1}, namely,
Theorem~6.6, which connects the principal specialization of the internal 
product of  Schur functions with certain parabolic Kostka polynomials, and
gives, via Corollary~6.7, \cite{Ki1},  an effective method for
computing the polynomials $s_{\al}*s_{\beta}(q,\ldots ,q^{N-1})$
which, turns out to be for the first time, does not use the character table 
of the  symmetric group $S_n$, ~$n=|\alpha|$. \\

Let $\alpha$ and $\beta$ be partitions,  $\ell(\alpha) =r,$ $\ell(\beta)=s$ 
~and $|\alpha|=|\beta|$. Let $N$ be an integer such that $r+s < N$. Consider 
partition
$$[\alpha, \beta]_N:=(\alpha_1+\beta_1,\alpha_2+\beta_1,\ldots,\alpha_r+\beta_1,
\underbrace{\beta_1,\ldots,\beta_1}_{N-r-s},\beta_1-\beta_s,\ldots,\beta_1-\beta_2].$$
Clearly, $| [\alpha,\beta]_N| = \beta_1 N$, $\ell([\alpha,\beta]_N=N-1$.
\begin{theorem}\label{t6.1} i) Let $\al ,\beta$ be partitions, $|\al|=|\beta|$,
$l(\al)\le r$, and $l(\al)+l(\beta)\le Nr$. Consider the sequence %$R_N$
of rectangular shape partitions
$$R_N=\{\underbrace{(\beta_1^r),\ldots ,(\beta_1^r)}_N\}.
$$
Then
\begin{equation}
s_{\al}*s_{\beta}(q,\ldots ,q^{N-1})\bdoteq K_{[\al
,\beta]_{{\bf Nr}},R_N}(q). \label{4.8}
\end{equation}

ii) (Dual form) Let $\alpha ,\beta$ be  partitions such that $|\alpha|=
|\beta |$, $ \al_1 \le r$ and $ \beta_1  \le k$. For given  integer $N$ such 
that $\alpha_1+\beta_1\le Nr$,  consider partition
$$\ld_N:=(rN-\beta_{k}^{'},rN-\beta_{k-1}^{'},\ldots ,rN-\beta_1^{'},\alpha^{'})
$$
and a sequence of rectangular shape partitions
$$R_N:=(\underbrace{(r^k),\ldots ,(r^k)}_N).$$
Then
\begin{equation}
K_{\ld_NR_N}(q)\bdoteq s_{\al}*s_{\beta}(q,\ldots ,q^{N-1})
\label{6.11*}
\end{equation}
\end{theorem}

\begin{theorem}~({\bf Fermionic formula for the principal specialization of the internal product of Schur functions}). ${}$

Let $\alpha$ and $\beta$ be two partitions of the same size, and 
$r:=\ell(\alpha)$ be the length of $\alpha$. Then

\begin{equation} 
s_{\alpha}*s_{\beta}(q,\ldots,q^{N-1})= \sum_{\{\nu \}} q^{c(\{\nu\})} \prod_{k,j \ge 1} { P_{j}^{(k)}(\nu)+m_j(\nu^{(k)})+N(k-1)\delta_{j,\beta_{1}} \theta(r-k) \brack P_j^{(k)}(\nu)}_{q},
\end{equation}
where the sum runs over the set of admissible configurations $\{\nu \}$ of 
type $([\alpha,\beta]_{N}, (\beta_1)^N)$. Here for any partition $\lambda$, 
~$\lambda_j$ denotes its $j$-th component.

\end{theorem}
Let us explain notations have used in Theorem~3.5.

$\bullet$~~A configuration $\{\nu \}$ of type $[\alpha,\beta]_N$ consists of 
a collection of partitions  $\{\nu^{(1)},\ldots,\nu^{(N-1)} \}$ such that 
$|\nu^{k}|= \sum_{j > k} \Bigl([\alpha,\beta]_N\Bigr)_j$; by definition we set
 $\nu^{(0)}:=(\beta_1)^N$; 

$\bullet$~~$P_{j}^{(k)}(\nu): = N \min(j,\beta_1)~ \delta_{k,1} + Q_j(\nu^{(k-1)}) -2 Q_j(\nu^{(k)}) +Q_j(\nu^{k+1})$; ~here for any partition $\lambda$ 
we set $Q_n(\lambda): =\sum_{ j \le n}~\min(n,\lambda_j)$;

$\bullet$~~For any partition $\lambda$,~~ $m_j(\lambda)$ denotes the number of 
parts of $\lambda$ are equal to $j$;  

$\bullet$~~A configuration $\{\nu\}$ of type $[\alpha,\beta]_N$ is called 
{\it admissible configuration of type $([\alpha,\beta]_N,(\beta_1)^{N})$}, if 
$P_{j}^{(k)}(\nu) \ge 0,~\forall j,k \ge 1 $;

$\bullet$~~Here ~$\delta_{n,m}$ denotes Kronecker's delta function, and we 
define $\theta(x) =1,~if~~ x \ge 0$, and $\theta(x)=0,~ if~~x <0$;

$\bullet$~~$c(\{\nu\})=\sum_{n,k \ge 1} {\lambda_n^{(k-1)}-\lambda_n^{(k)}
 \choose 2}$ denotes the {\bf charge} 
of a configuration $\{ \nu \}$; by definition,  ${x \choose 2}:=x(x-1)/2,~\forall x \in \R$. \\

Let us draw attention to the fact that the summation in $(4.9)$ runs over the 
set of all admissible configurations of type $([\alpha,\beta]_N,(\beta_1)^N)$,
 other than that of type $([\alpha,\beta]_{Nr},(\beta_{1}^{r})^N)$.

\begin{cor} (\cite{Ki2},\cite{Ki3})${}$

For any partitions of the same size $\alpha$ and $\beta$, the polynomial 
$s_{\alpha}*s_{\beta}(q,\ldots,q^{N-1})$ is symmetric and {\bf unimodal}. 
In particular, the 
 generalized Gaussian polynomial ${N \brack \alpha}_q$ is symmetric and 
{\bf unimodal} for any partition $\alpha$.
\end{cor}
Indeed, in the case $\beta=(n)$,  $n:= |\beta|$, one has 
$s_{\alpha}*s_{(n)}=s_{\alpha}$.

Our proof of Corollary~3.7 is based on the following identity
$$2c(\nu)+\sum_{k,j \ge 1} P_{j}^{(k)}(\nu)\biggl[m_j(\nu^{(k)})+N(k-1) 
\delta_{j,\beta_{1}} \theta(r-k) \biggr] = N |\alpha|,$$
which can be checked directly. This identity shows that the all polynomials 
associated with a given admissible configuration involved, are symmetric and 
have the same ``symmetry center'' $N |\alpha|/2$, and therefore the resulting 
polynomial $s_{\alpha}*s_{\beta}(q, \ldots,q^{N-1})$ is symmetric and \underline{unimodal}. 
 \begin{cor} Let $\alpha$ and $\beta$  be partitions of the same size, and 
$K_{\beta,\lambda}(q,t)$ denotes the Kostka--Macdonald polynomial associated 
with partitions $\alpha$ and $\beta$, \cite{Ma}. One has
$$ K_{\beta,\alpha}(q,q)= H_{\alpha}(q) \Biggl( \sum_{\{\nu\}} q^{c(\{\nu\})}~
\prod_{k=2}^{r} \frac{1}{\biggl[m_{\beta_{1}}(\nu^{(k)})\biggr]_{q} !} ~~\prod_{k \ge 1 \atop j \ge 1, ~\j \not= \beta_{1}} {P_{j}^{(k)}(\nu)+m_{j}(\nu^{(k)}) \brack m_{j}(\nu^{(k)})}_{q} \Biggr), $$ 
where the sum runs over the same set of admissible configurations as in Theorem~3.6, and $[m]_{q} ! := \prod_{j=1}^{m} (1-q^{j})$ stands for the $q$-factorial
 of an positive integer $m$, and by definition $[0]_{q} !=1$; $H_{\alpha}(q):=
\prod_{x \in \alpha}(1-q^{h(x)})$ denotes the hook polynomial associated with partition $\alpha$, see e.g. \cite{Ma}.
\end{cor}
Indeed, one can show \cite{Ki4}, \cite{Ki2} that
$$ \lim_{N \rightarrow \infty} s_{\alpha}*s_{\beta} (1,q,q^2,,q^{N}) =
\frac{K_{\beta,\alpha}(q,q)}{H_{\alpha}(q)} = \lim_{N \rightarrow \infty}K_{[\alpha,\beta]_N, (\beta_{1})^{N}}(q).$$

A fermionic formula for the principal specialization of the internal product 
of Schur functions, and therefore that for the generalized Gaussian polynomials, is a far generalization of the so-called $KOH$-identity \cite{O} which is equivalent to the fermionic formula for the Kostka number $K_{(Nk,k),(k)^{N+1}})(1)$. The rigged configuration bijection gives rise to a combinatorial proof of 
Theorem~3.5, and therefore to a combinatorial proof of  unimodality of the 
generalized Gaussian polynomials \cite{Ki2}, as well as to give an 
interpretation of the statistics introduced in \cite{O} in terms of rigged 
configurations data, see \cite{Ki4}, Section~10.2.

\begin{ex}  Let $\alpha =(4,2)$,~ $\beta=\alpha ^{'} =(2,2,1,1).$ We want to 
compute the principal specialization of the internal product of Schur functions ~~$s_{\alpha}*s_{\beta}(q,\ldots,q^{N-1})$ by means of a fermionic formula 
$(3.10)$.First of all, there are $8$ admissible configurations of type 
$((4,2),(2,2,1,1),(6)^N)$.In fact, it is a general fact that for given partitions $\lambda$ and $\mu$,  the number of admissible configurations of type $(N \lambda,N \mu)$ {\bf doesn't depend} on $N$ ,if $N > N_{0}$ i a certain  number 
$N_0:=N_{0}(\lambda,\mu)$ depending on $\lambda$ and $\mu$ only.  This fact is 
a direct consequence of constraints are imposed by the  set of inequalities 
$\{ P_j^{(k)}(\nu) \ge 0,~~\forall j,k \ge 1 \}$. \\ 
Now let us list the {\bf conjugate} of the \underline{first configurations} 
$\nu^{(1)} \in \{\nu \}$  for all admissible configurations  $\{\nu\}$ of type $[(4,2),(2,2,1,1)_N]$,
together with all non-zero numbers $P_{j}^{(k)}(\nu),~~j,k \ge 1$.
$$(N-3,N-3),~~P_{2}^{(1)}=2,P_{2}^{(2)}=2,~~c=9,$$
$$(N-2,N-4),~~P_{2}^{(1)}=2,P_{1}^{(2)}=1, P_{2}^{(2)}=2, ~~9,$$            
$$(N-3,N-4,1),~~P_{1}^{(1)}=2,P_{1}^{(2)}=4,P_{1}^{(3)}=2, P_{2}^{(2)}=1,~~ c=11,$$
$$(N-2,N-5,1),~~P_{1}^{(2)}=1,P_{2}^{(1)}=4,P_{3}^{(1)}=2, P_{2}^{(2)}=1, ~~ c=13,$$
$$(N-3,N-5,1,1),~~P_{1}^{(1)}=2,P_{2}^{(1)}=6,(P_{2}^{(2)}=0),P_{3}^{(1)}=2, ~~ c=15,$$
$$(N-3,N-5,2),~~P_{1}^{(1)}=2,P_{2}^{(1)}=6,(P_{2}^{(2)}=0),P_{3}^{(1)}=2, ~~ c=17,$$
$$(N-2,N-6,1,1),~~P_{1}^{(2)}=1, P_{2}^{(1)}=6,(P_{2}^{(2)}=0),P_{4}^{(1)}=2, ~~ c=19,$$
$$(N-2,N-6,2),~~P_{1}^{(2)}=1, P_{2}^{(1)}=6,(P_{2}^{(2)}=0), P_{3}^{(1)}=2, ~~ c=21.$$
These are data related with the first configurations $\nu^{(1)}$ in the set 
of all admissible configurations of type $([(42),(2211)]_{N},(6)^{N})$, and 
the set of  all nonzero numbers $P_{j}^{(k)}(\nu)$.~~ All 
other diagrams~ $\nu^{(k)},~k > 1$~ from the set of admissible configurations 
in question, are the \underline{same}, and are displayed below
$$\nu^{(k)}=(N-2-k,max(N-4-k,0)),~~2 \le k \le N-3,$$
so that ~~~$m_{1}^{(2)}=2,~~m_{2}^{(2)} =N-6$.

Therefore,
$$ (\clubsuit)~~~K_{[(4,2),(2,2,1,1)]_{2 N},(2,2)^N}(q) \bdoteq  q^9 {N-1 \brack 2}~{2 N-4 \brack 2}+
q^9 {3 \brack 1} {N-2 \brack 2} {2 N-4 \brack 2} +$$
$$q^{11} {3 \brack 1} {3 \brack 1} {N-1 \brack 4} {2 N-5 \brack 1} + 
q^{13} {3 \brack 1} {3 \brack 1} {N-2 \brack 4} {2 N-5 \brack 1} +$$
$$ q^{15} {3 \brack 1} {4 \brack 2} {N \brack 6}+q^{17} {4 \brack 2} 
{4 \brack 2} {N-1 \brack 6}+ q^{19} {3 \brack 1} {3 \brack 1} {N-1 \brack 6} 
+ q^{21} {3 \brack 1} {4 \brack 2} {N-2 \brack 6}.$$
On the other hand,
$$s_{42}*s_{2211}=\frac{1}{720} \Biggl( 81 p_{1}^6 - 135 p_{1}^4 p_{2} +
45 p_{1} p_{2}^2 -90 p_{1}^2 p_{4}+144 p_{1} p_{5}-135 p_{2}^3 +90 p_{2} p_{4} \Biggr).$$ 
Here ~$p_{k}:= \sum_{i \ge 1} x_{i}^{k}$ stands for the power sum symmetric 
functions degree of $k$. One can check that $s_{42}*s_{2211}(q,\ldots,q^{N-1})
 =  RHS(\clubsuit) \bdoteq K_{[(4,2),(2,2,1,1)]_{2 N},{(2,2)}^N}(q)$, as 
expected. \\
Finally one can check that ~~$\lim_{N \rightarrow \infty}~s_{42}*s_{2211}(q,
\ldots,q^{N-1}) =q^9(2,1,2,2,1,1) =K_{2211,42}(q,q)$.\\
\end{ex}

\subsection{Polynomiality of stretched Kostka and Littlewood--Richardson numbers} 

 %{ \rm Few comments in order.
{\rm
 $\bullet$~~ As it was mentioned above, for a given partitions $\lambda$ and 
$\mu$ (resp. $\lambda$ and a sequence of rectangular shape partitions 
$\{\cal{R}\}$), the  number of {\bf admissible} configurations of type 
$(N \lambda, N \mu)$ (resp. of type  $(N \lambda, \{\cal{R}\}$))  
{\bf doesn't}  depend on $N > N_{0}$, where the number  $N_{0}$ depends on 
$\lambda$ and $\mu$ (resp. $\lambda$ and $\{\cal{R}\}$) only.

$\bullet$~~Each admissible configuration $\{\nu\}$ provides a contribution of 
a form 
$$ a_{{\{\nu}\}}(q) \prod_{j,k\ge 1} {b_{j,k}(\nu) N + d_{jk}(\nu) \brack 
d_{jk}(\nu)}_{q}, $$
to the parabolic Kostka polynomial $K_{N \lambda, N \{\cal{R}\}}(q)$,~~where 
a polynomial $a_{\nu}(q)$ and a {\bf finite set} of numbers 
$\{b_{jk}(\nu), d_{jk}(\nu) \}_{ j,k \ge 1}$ both doesn't depend on $N > N_{0}$
for some $N_{0}:=N_{0}(\lambda,\{ \cal{R} \})$.

$\bullet$ It is clear that the sum~~$\sum_{N \ge 0} {a N+b \choose b} t^{N}$
~~ is a  {\bf rational} function of variable $t$. It is well-known (and easy 
to prove)  that the {\it Adamar product}
\footnote{~~ See e.g. $wikipedia.org/wiki/Hadamard{\_}product{\_}(matrices)$ 
and the literature quoted therein.}
of {\it rational} functions is again a {\it rational} function. Therefore,
\begin{cor} ~(\cite{Ki1},\cite{KT}, \cite{Ras}) ${}$

For any partition $\lambda$ and a sequence of {\it rectangular shape 
partitions} $\{\cal{R}\}$,  the generating function 
$$\sum_{N \ge 0}~K_{\lambda,\{\cal{R}\}}(1)~ t^N$$
is a rational function of variable $t$ with a unique pole at~~ $t=1$. 
\end{cor}
More generally using a $q$-version of Adamar's product Theorem, we can show
\begin{theorem} (\cite{Ki1}) ${}$

For any partition $\lambda$ and a \underline{dominant} sequence of 
 rectangular shape partitions ~ $\{\cal{R}\}$, ~ the generating function of 
stretched parabolic Kostka polynomials 
$$\sum_{N \ge 0}~K_{\lambda,\{\cal{R}\}}(q)~ t^N$$
is a rational function of variables $q$ and $t$ of a form 
$P_{\lambda,\{\cal{R\}}}(q,t) / Q_{\lambda,\{\cal{R}\}}(q,t)$,~where the 
dominator $ Q_{\lambda,\{\cal{R}\}}(q,t)$ has the following form
$$Q_{\lambda,\{\cal{R}\}}(q,t) = \prod_{ s \in S}(1-q^s~t)$$
for a certain finite set $S:=S(\lambda, \{\cal{R}\})$ depending on $\lambda$ 
and  $\{\cal{R}\}$. 
\end{theorem}
Clearly that Corollary~3.9 is a special case $q=1$ of Theorem~3.10.

$\bullet$ ~~(Littlewood--Richardson polynomials)~~Let $\lambda$ be a partition and $\{\cal{R}\}$ be a dominant sequence of rectangular shape partitions. 
~Write
$$ K_{\lambda,\{\cal{R}\}}(q)= {b(\lambda, \cal{R})}~q^{a(\lambda, \cal{R})} +~~higher~~degree~~terms.$$

$(1)$~~(Generalized saturation theorem \cite{Ki3})
$${a(N \lambda, N \{\cal{R}\})} = N~ a(\lambda,  \{\cal{R}\}).$$ 
Therefore,
$$ \sum_{N \ge 0} {b(N \lambda, N  \{\cal{R} \})}~t^N = \frac{ {P_{\lambda , \cal{R}}}(q, q^{-a(\lambda,  \{\cal{R}\})}~t)}{ {Q_{\lambda,\cal{R}}}(q,q^{-a(\lambda,  \{\cal{R}\})}~ t)} \Bigg\vert_{q=0} $$
is a {\bf rational} function with a unique pole at $t=1$ with multiplicity 
equals to $\#| s \in {S(\lambda, \cal{R})} ~|~ s=a(\lambda,{\cal{R}}) |$.

$(2)$~~~Now let $\lambda$,~$\mu$ ~and~$\nu$ be partitions such that $|\lambda|+|\mu|=|\nu|$. Consider an integer $N \ge max(\ell(\lambda), \mu_1)$, and 
define  partition $\Lambda = \Lambda(N,\lambda,\mu):= (N^N) \oplus \lambda, 
\mu)$ and the dominant rearrangement of the set of rectangular shape 
partitions  $\{ (N^N),\nu_1,\ldots, \nu_{\ell(\nu)} \}$, denoted by $M:= 
M(N,\nu)$.
\begin{pr} (\cite{Ki3}) ${}$ 

One has
$$ b(\Lambda, M) := c_{\lambda, \mu}^{\nu},$$
where $c_{\lambda, \mu}^{\nu}$ denotes the Littlewood--Richardson number 
corresponding to partitions $\lambda, \mu$~and~$\nu$, that is , the 
multiplicity of Schur function $s_{\nu}$ in the product of Schur functions 
$s_{\lambda}~s_{\mu}$.
\end{pr}
\begin{theorem}(\cite{Ki3},\cite{Ras}) ${}$

Given three partitions $\lambda, \mu$~and~$\nu$ such that $|\lambda|+|\mu|=
|\nu|$.
~~The the generating function
$$ \sum_{N \ge 0}~c_{N \lambda,N \mu}^{N \nu}~t^N $$
is a {\bf rational} function of variable $t$ with a unique pole at $t=1$. 
Therefore, $c_{N \lambda,N \mu}^{N \nu}$ is a {\bf polynomial} in $N$ with 
rational coefficients.
\end{theorem}
\begin{ex} (\cite{Ki4}, \cite{Ki1})~~(MacMahon polytope and multidimensional 
Narayana numbers again) ${}$

Take $\ld =(n+k,n,n-1,\ldots ,2)$ and $\mu
=\ld'=(n,n,n-1,n-2,\ldots ,2,1^k)$. One can show \cite{Ki1}  that if 
$n\ge k\ge 1$, then for any positive integer $N$
\vskip 0.2cm

$\bullet$ $a(N\ld,N\mu)=(2k-1)N$;

$\bullet$ $b(N\ld,N\mu)={\rm dim}V_{((n-k+1)^{k-1})}^{\g
l(N+k-1)}=\ds\prod_{i=1}^{k-1}\prod_{j=1}^{n-k+1}\frac{N+i+j-1}{i+j-1}.$

In other words, the number $b(N\ld,N\mu)$ is equal to the number of (weak)
plane partitions of rectangular shape $((n-k+1)^{k-1})$ whose
parts do not exceed $N$. According to Exercise~1, {\bf c}, \cite{Ki1}, pp. 102--103, $b(N\ld,N\mu)$ is equal also to the number
$i(\M_{k-1,n-k+1};N)$ of rational points ${\bf x}$ in the
MacMahon polytope $\M_{k-1,n-k+1}$ such that the points $N{\bf x}$
have integer coordinates. 
%Therefore, in this example one can take
%$$\Gamma (\ld,\mu)=\M_{k-1,n-k+1}.$$
It follows from (2.8) that the generating function for
numbers $b(n\ld,n\mu)$ has the following form
$$\sum_{n\ge 0}b(n\ld,n\mu)t^n=\left(\sum_{j=0}^{(k-2)(n-k)}
N(k-1,n-k+1;j)t^j\right)/(1-t)^{(k-1)(n-k+1)+1},
$$
where $N(k,n;j)$, $0\le j\le (k-1)(n-1)$, denote   rectangular
Narayana's numbers, see e.g. \cite{Mc}, \cite{Su1}.

One can show (A.K.)  that

$\bullet$ ~  if $r:=k- {n+2 \choose 2} \ge 0$, then $b(\ld,\mu)=1$, and
$$a(\ld,\mu)=2 { n+3 \choose 3}+(n+1)(2 r-1)+ {r \choose 2};$$

$\bullet$ if $1\le k< {n+2 \choose 2}$, then there exists a unique $p$, 
$1\le p\le n$, such that
$$(p-1)(2n-p+4)/2<k\le p(2n-p+3)/2.$$
In this case
$$a(\ld,\mu)=p(2k-(p-1)n-p)+2 {p \choose 3}, $$
and one can take $\Gamma(\ld,\mu)$ to be
equal to the MacMahon polytope $\M_{r(k),s(k)}$ with \\
$r(k):=k-1-(p-1)(2n-p+4)/2$, and $s(k):=p(2n-p+3)/2-k$.

This Example  gives some flavor how intricate the piecewise
linear function $a(\ld,\mu)$ may be.
\begin{con}

Let $\lambda$ and $\mu$ be partitions of the same size. Then

$\bullet$~~(\cite{Ma})~~~~~~~$a(\lambda,\mu) = a(\mu^{'}, \lambda^{'})$,

$\bullet$~~(\cite{Ki4})~~~~~~$b(\lambda,\mu) = b(\mu^{'}, \lambda^{'})$.
\end{con}

\qed

\end{ex}

\begin{de} (\cite{Ki1})~~~ Let $\alpha$ and $\beta$ be partitions of the same 
size. Define 
Liskova polynomials $L_{\alpha, \beta}^{\mu}(q)$ through the decomposition of 
the internal product of Schur functions in terms of Hall-Littlewood polynomials
$$ s_{\alpha}*s_{\beta}(X) = \sum_{\mu} L_{\alpha, \beta}^{\mu}~P_{\mu}(X;q).$$
Clearly, $L_{\alpha,\beta}^{\mu} \in \N[q]$, and $L_{\alpha,(|\alpha||}(q)=
K_{\alpha,\mu}(q)$, so that the Liskova polynomials are natural generalization 
of Kostka--Foulkes polynomials.
\end{de}
\begin{prb} Find for Liskova polynomials an analogue of a fermionic formula 
for Kostka--Foulkes polynomials stated, for example, in \cite{Kir}, \cite{Ki4}.      
\end{prb}
}

\section{ Appendix. {\bf  Rigged Configurations}: a brief review} 

Let $\ld$ be a partition and $R=((\mu_a^{\eta_a}))_{a=1}^p$ be a sequence
of rectangular shape partitions such that
$$|\ld|=\sum_a|R_a|=\sum_a\mu_a\eta_a.$$

\begin{de}  ${}$

The configuration of type $(\ld ,R)$ is a sequence
of {\bf partitions}  $\{\nu \}=(\nu^{(1)},\nu^{(2)},\ldots )$ such that
$$|\nu^{(k)}|=\sum_{j>k}\ld_j-\sum_{a\ge 1}\mu_a\max (\eta_a-k,0)
=-\sum_{j\le k}\ld_j+\sum_{a\ge 1}\mu_a\min (k,\eta_a)
$$
for each $k\ge 1$.
\end{de}

Note that if $k\ge l(\ld)$ and $k\ge\eta_a$ for all
$a$, then $\nu^{(k)}$ is empty.

 In the sequel we make the convention that  $\nu^{(0)}$ is the empty 
partition
\footnote{~~ However, in some cases it is more convenient to set $\nu^{(0)}=
(\mu_{i_{1}},\ldots,\mu_{i_{s}})$, where we assume that $\eta_{i_{a}}=1, ~a=1,
\ldots,s$. We will give an indication of such choice if it is necessary.
}.\\
For a partition $\mu$ and an integer $j \ge 1$  define the number 
$$Q_j(\mu)=\mu_1'+\cdots +\mu_j',$$
 which is equal to the number of cells in the first $j$  columns of $\mu$. 

The {\it vacancy} numbers ~$P_j^{(k)}(\nu
):=P_j^{(k)}(\nu ;R)$ ~ of the configuration $\{\nu \}$ of type $(\ld,R)$ are 
defined by
$$P_j^{(k)}(\nu )=Q_j(\nu^{(k-1)})-2Q_j(\nu^{(k)})+Q_j(\nu^{(k+1)})
+\sum_{a\ge 1}\min (\mu_a,j)\delta_{\eta_a,k}
$$
for $k,j\ge 1$, where $\delta_{a,b}$ is the Kronecker delta.

\begin{de}\label{d4.2} The configuration $\{\nu \}$ of type $(\ld ,R)$ is 
called  {\bf admissible}, if
$$P_j^{(k)}(\nu ;R)\ge 0~~{\rm for~all}~~ k,j\ge 1.$$
\end{de}

We denote by $C(\ld ;R)$ the set of all admissible configurations
of type $(\ld ,R)$, and call the vacancy number $P_j^{(k)}(\nu,R)$ 
{\it essential}, if $m_j(\nu^{(k)})>0$.

Finally, for configuration $\{\nu \}$ of type $(\ld ,R)$ let us define its
{\bf charge}
$$ c(\nu )=\sum_{k,j\ge 1}\left(\begin{array}{c}\alpha_j^{(k-1)}-
\alpha_j^{(k)}+\sum_a\theta (\eta_a-k)\theta (\mu_a-j)\\ 2\end{array}
\right),
$$
and {cocharge}
$$\overline c(\nu )=\sum_{k,j\ge 1}\left(\begin{array}{c}\alpha_j^{(k-1)}
-\alpha_j^{(k)}\\ 2\end{array}\right),
$$
where $\alpha_j^{(k)}=(\nu^{(k)})_j'$ denotes the size of the $j$--th
column of the $k$--th partition $\nu^{(k)}$ of the configuration $\{\nu \}$;
for any real number $x\in \R$ we put $\theta (x)=1$, if $x\ge 0$, and
$\theta (x)=0$, if $x<0$.

\begin{theorem}\label{t4.3} {\rm ({\bf Fermionic formula for parabolic Kostka
polynomials} \cite{Ki4})}  ${}$

Let $\ld$ be a partition and $R$
be a dominant sequence of rectangular shape partitions. Then
\begin{equation}
K_{\ld R}(q)=\sum_{\nu}q^{c(\nu )}\prod_{k,j\ge 1}\left[\begin{array}{c}
P_j^{(k)}(\nu ;R)+m_j(\nu^{(k)})\\ m_j(\nu^{(k)})\end{array}\right]_q,
\label{4.1}
\end{equation}
summed over all admissible configurations $\nu$ of type $(\ld ;R)$;
$m_j(\ld)$ denotes the number of parts of the partition $\ld$ of size
$j$.
\end{theorem}

\begin{cor}\label{c4.4} {\rm ({\bf Fermionic formula for Kostka--Foulkes
polynomials \cite{Kir}})}  ${}$

Let $\ld$ and $\mu$ be partitions of the same
size. Then
\begin{equation}
K_{\ld\mu}(q)=\sum_{\nu}q^{c(\nu)}\prod_{k,j\ge 1}\left[\begin{array}{c}
P_j^{(k)}(\nu ,\mu)+m_j(\nu^{(k)})\\ m_j(\nu^{(k)})\end{array}\right]_q,
\label{4.1a}
\end{equation}
summed over all sequences of partitions $\nu
=\{\nu^{(1)},\nu^{(2)},\ldots\}$ such that

$\bullet$ ~~$|\nu^{(k)}|=\sum_{j>k}\ld_j$, $k=1,2,\ldots$;

$\bullet$ ~~ $P_j^{(k)}(\nu ,\mu):=Q_j(\nu^{(k-1)})-2Q_j(\nu^{(k)})
+Q_j(\nu^{(k+1)})\ge 0$  for all $k,j\ge 1$, ~~where by definition we
\underline{put} $\nu^{(0)}=\mu$;

\begin{equation}
\bullet~~c(\nu )=\sum_{k,j\ge
1}\left(\begin{array}{c}(\nu^{(k-1)})_j' -(\nu^{(k)})_j'\\
2\end{array}\right).~~~~~~~~~~~~~~~~~~~~~~~~~~~~~\label{4.3*}
\end{equation}
\end{cor}

It is frequently convenient to represent an admissible configuration $\{\nu\}$
by a matrix $m(\nu)=(m_{ij}),~m_{ij} \in \Z, \forall i,j \ge 1$, which must 
meets certain conditions. Namely, 
starting from the collection of partitions $\{\nu \} =(\nu^{(1)}, \nu^{(2)},
 \ldots, ...)$ corresponding to configuration $\{\nu\}$,  define matrix
$$ m(\nu):= (m_{ij}),~~m_{ij}=(\nu^{(i-1)})_{j}^{'} -(\nu^{(i)})_j^{'}+ 
\sum_{a \ge 1} \theta (\eta_{a}-i) \theta (\mu_{a}-j),~~\nu^{(0)}:= 
\emptyset, $$
where we set by definition ~$\theta(x)= 1,~if ~x \in \R_{\ge 0}$ and $\theta(x)=0,~x \in \R_{< 0}$. \\
One can check that a configuration $\{\nu\}$ of type $(\ld, R)$ is 
{\bf admissible} if and only if the matrix $m(\nu)$ meets the following conditions

$(0)$~~$m_{ij} \in \Z$,

$(1)$~~$\sum_{i \ge 1}m_{ij} = \sum_{a \ge 1} \eta_{a} \theta (\mu_{a}-j),$

$(2)$~~$\sum_{j \ge 1} m_{ij} = \ld_{i}$, 

$(3)$~~$\sum_{j \le k} (m_{ij} -m_{i+1,j}) \ge 0$, for all $i,j,k$

$(4)$~~$\sum_{a \ge 1} \min(\eta_{a},k) \delta_{\mu_{a}, j} \ge \sum_{i \le k}
(m_{ij}- m_{i,j+1})$, for all $i,j,k$.

One can check that if matrix $(m_{ij})$  satisfies the conditions $(0)-(4)$, 
then the set of partitions  $\{\nu\} = (\nu^{(1)}, \nu_{(2)}, \ldots,...)$ 
$$ (\nu^{(k)})_{j}^{'}:=\sum_{i > k} m_{ij} - \sum_{a} \max(\eta_{a}-k,0) 
\theta(\mu_{a}-j)$$
defines an admissible configuration of type $(\ld,R=\{(\mu_{a})^{\eta_{a}}\})$.

\begin{ex} Take $\ld =(44332)$, ~$R = \{ (2^3),(2^2),(2^2),(1),(1) \}$, so that
$$\{\mu_a\}=(2,2,2,1,1)~~ and~~ \{\eta_a \}=(3,2,2,1,1),~~a=1,\ldots, 5.$$.
Therefore 
$|\nu^{(1)}|=4$, $|\nu^{(2)}|=6$, $|\nu^{(3)}|=5$, and $|\nu^{(4)}|=2$. It is 
not hard to check that there exist $6$ admissible configurations.  They are: \\
$(1)$~ $\{\nu^{(1)}= (3,1),\nu^{(2)}=(3,3),\nu^{(3)}=(3,2),\nu_{(4)}=(2) \}$,\\
$(2)$~ $\{\nu^{(1)}= (3,1),\nu^{(2)}=(3,2,1),\nu^{(3)}=(3,2),\nu^{(4)}=(2) \}$,\\
$(3)$~ $\{\nu^{(1)}= (2,2),\nu^{(2)}=(2,2,2),\nu^{(3)}=(3,2),\nu^{(4)}=(2) \}$,\\
$(4)$~ $\{\nu^{(1)}= (4),~\nu^{(2)}=(3,3),\nu^{(3)}=(3,2),\nu^{(4)}=(2) \}$, \\
$(5)$~ $\{\nu^{(1)}= (3,1),\nu^{(2)}=(2,2,1,1),\nu^{(3)}=(2,2,1),\nu^{(4)}=(2) \}$, \\
$(6)$~ $\{\nu^{(1)}= (3,1),\nu^{(2)}=(2,2,1,1),\nu^{(3)}=(3,1,1),\nu^{(4)}=(2) \}$, \\
Let us compute the matrix $(m_{ij})$ corresponding to the configuration $(2)$. 
Clearly, \\
$(m_{ij}) =((\nu^{(i-1)})_j^{'} - (\nu^{(i)})_j^{'}) +(\sum_{a \ge 1}
\theta(\eta_a-i) \theta(\mu_a-j)) := U +W$.
%\end{ex}
One can check that
$$ U=  \left( \begin{array}{ccccc} 
-3 & -1 & 0 & 0 & 0 \\
0 & -1 &-1 & 0 & 0\\
0 & 0 & 1 & 0 & 0\\
2 & 1 & 0 & 0 & 0\\
1 & 1 & 0 & 0 & 0
\end{array} \right), ~~~
 W= \left( \begin{array}{ccc}
5 & 3 & 0  \\
3 & 3 & 0 \\
1 & 1 & 0
\end{array} \right).
$$
Therefore,
$$ m(\{ \nu \}) = \left( \begin{array}{ccccc}
2 & 2 & 0 & 0& 0\\
3 & 2 & -1 & 0 & 0 \\
1 & 1 & 1 & 0 & 0 \\
2 & 1 & 0 & 0 & 0 \\
1 & 1 & 0 & 0 & 0
\end{array} \right).
$$
{ \rm One can read off directly from the matrix $m_{ij}$ the all  additional 
quantities need to compute the parabolic Kostka 
polynomial corresponding to $\lambda$ and a (dominant) sequence of 
rectangular shape partitions $R$. Namely,
$$P_{j}^{(k)}= \sum_{i \ge j}(m_{ki}-m_{k+1,i}),~~m_{j}(\nu^{(k)}) = 
\sum_{a \ge 1} \min(\eta_{a},k) \delta_{\mu_{a}, j} -  \sum_{i \le k}
(m_{ij}- m_{i,j+1}), c(\nu)=\sum_{i,j \ge 1} {m_{ij} \choose 2}.$$
For example, in our example, we have $c(\nu)=8,$~$P_{1}^{(1)}=1,P_{2}^{(2)}=1,
P_{3}^{(2)}=1,P_{2}^{(3)}=1$ are all non-zero vacancy numbers, and the contribution of the configuration in question to the parabolic Kostka polynomial is equal to 
$q^8 {2 \brack 1}^{4}$. Treating in a similar fashion other configurations, we come to a fermionic formula
$$K_{44332,\{(2^3),(2^2),(2^2),(1),(1)\}}(q)= q^{10} {3 \brack 1}+q^8 {2 \brack 1}^4+q^{8} {3 \brack 2}+q^{12}+q^6 {2 \brack 1}{3 \brack 2} + q^8.$$
}
\end{ex}

\qed

If $\eta_{a}=1,~\forall a$, then  $\sum_{a \ge 1} \eta_{a} \theta(\mu_a -j) = 
\mu_{j}^{'}$, and we \underline{set} $\nu^{(0)}=\mu^{'}$. In this case one can rewrite the  conditions $(1)-(4)$ as follows

$(1^{'})$~~$\sum_{i \ge 1} m_{ij} = \mu_{j}^{'}$,

$(2^{'})$~~$\sum_{j \ge 1} m_{ij} =\ld_{i} $, 

$(3^{'})$~~$\sum_{j \le k} (m_{ij} -m_{i+1,j}) \ge 0$, for all $i,j,k$

$(4^{'})$~~$\sum_{i > k} (m_{ij} -m_{i,j+1}) \ge 0$, for all $i,j.k$

Let us remark that if $m_{ij} \in \Z_{\ge 0}$ then the matrix $(m_{ij})$ defines a 
lattice plane partition of shape $\ld$.
 For example, take $\ld=(6,4,2,2,1,1)$, $\mu=(2^8)$ and admissible configuration 
$\{ \nu\}=\{ (5,5),(4,2),(3,1),(2),(1) \}$. The corresponding matrix and lattice 
plane partition of shape $\ld$ are
$$ (m_{ij}) =\left( \begin{array}{cccccc}
3 & 3 & 0 & 0 & 0 & 0 \\
1 & 3 & 0 & 0 & 0 & 0 \\ 
1 & 1 & 0 & 0 & 0 & 0 \\
1 & 1 & 0 & 0 & 0 & 0 \\
1 & 0 & 0 & 0 & 0 & 0 \\ 
1 & 0 & 0 & 0 & 0 & 0 \end{array} \right),~~~~~~~
and~~plane~~partition~~\begin{array}{cccccc}
3 & 3 \\
1 & 3 \\
1 & 1 \\
1 & 1 \\
1 \\
1 \end{array}.
$$ 
The corresponding lattice word is $111222.1222.12.12.1.1$. \\
  
In the case $\eta_{a}=1,~\forall a$, there exists a  unique admissible 
configuration of type $(\lambda,\mu)$ denoted by $\Delta(\ld,\mu)$ such that 
$\max( c((\Delta(\ld,\mu), J)))= n(\mu) - n(\ld)$, where the maximum is taken 
over all     \underline{rigged configurations}  associated with configuration 
$\Delta(\ld,\mu)$. Recall that for any partition $\ld$, ~
$$n(\lambda) = \sum_{j \ge 1} {\ld_{a}^{'} \choose 2},$$
and if $\ld \ge \mu$ with respect to the {\it dominance order}, then the degree of 
Kostka polynomial $K_{\lambda,\mu}(q)$ is equal to ~~$n(\mu)-n(\ld)$,~~ see, e.g.  
\cite{Ma}, Chapter 1, for details. Namely, the configuration $\Delta(\ld,\mu)$ 
corresponds to the following matrix:
$$m_{1j}= \mu_{j}^{'} -\max(\ld_{j}^{'} -1,0), j \ge 1,~~~m_{ij}=1, ~~if~~(i,j) \in \ld, ~i \ge 2,~~m_{ij}=0,~if~(i,j) \notin \ld.$$
In other words, the configuration $\Delta(\ld, \mu)$ consists of the following 
partitions ~$(\ld[1],\ld[2], \ldots)$, where $\ld[k] = (\ld_{k+1},\ld_{k+2},\ldots)$. It is not difficult to see that the contribution to the Kostka polynomial $K_{\ld,\mu}(q)$ coming from the maximal configuration, is equal to
$$ K_{q}(\Delta(\ld, \mu)):=  q^{c(\Delta(\ld,\mu))} ~\prod_{j= 1}^{\ld_{2}} {Q_{j}(\mu) -Q_j(\ld)+\ld_{j}^{'}- \ld_{j+1}^{'} \brack \ld_{j}^{'}- \ld_{j+1}^{'} }_{q},$$
where $c(\Delta(\ld,\mu)) = n(\ld)+n(\mu)-\sum_{j \ge 1} ~\mu_{J}^{'} (\ld_{j}^{'} -1)$.~~  Therefore,
\begin{equation} 
K_{\ld, \mu}(q) \ge K_{q}(\Delta(\ld,\mu)).
\end{equation}
It is clearly seen that if $\ld \ge \mu$, then $Q_{j}(\mu) \ge Q_{j}(\ld),~\forall~ j \ge 1$, and thus, $K_{q=1}(\Delta(\ld,\mu)) \ge 1$, and the inequality 
$(5.13)$ can be considered as a ``quantitative'' generalization of the Gale--
Ryser theorem, see, e.g. \cite{Ma}, Chapter~ I, Section~7, or \cite{Ki8} for 
details. \\

Now let us {\bf stress} that for a fixed $k$, the all partitions $\nu^{(k)}$ 
which contribute to the set of admissible configurations of type $(\ld, \mu)$ 
have the same size equals to $\sum_{j \ge k+1} \ld_j$, and thus the size of 
each $\nu^{(k)}$  doesn't depend on $\mu$. However the Rigged Configuration 
bijection 
$$RC_{\ld,\mu}:  STY( \ld, \mu) \longrightarrow RC(\ld, \mu)$$
happens to be essentially  depends on $\mu$. One can check that the map 
$RC_{\ld,\mu}$ is compatible with the familiar {\it Bender--Knuth} 
transformations on the set of semistandard Young tableaux of a fixed shape.

As it was mentioned above, for a fixed $k$ the all (admissible) configurations
 have the same size. Therefore, the set of admissible configurations admits 
a partial ordering  denoted by ~``$ \succcurlyeq$'' .  ~~Namely, if $\{\nu\}$  and 
 $\{\xi\}$ are two admissible configurations of the same type $(\ld,\mu)$, 
we will write $\{\nu \} \ \succcurlyeq \{\xi\}$, if either $\{\nu\} =\{\xi\}$ 
or there exists an integer $\ell$ such that $\nu^{(a)}= \xi^{(a)}$~if $1 \le 
a \le \ell$, and   $\nu^{(\ell+1)} > \xi^{(\ell+1)}$ with respect of the 
dominance order on the set of the same size partitions. It seems an 
interesting 
{\bf Problem} to study  poset structureds on the set of  admissible 
configurations of type $(\ld,\mu)$, especiallyto investigate  the posets of 
admissible configurations associated with the multidimensional Catalan numbers, (work in progress).
 
\qed

\begin{theorem} ({\bf Duality theorem for parabolic Kostka polynomials} 
\cite{Ki4}) ${}$

Let $\lambda$ be partition and $R=\{ (\mu_{a}^{\eta_{a}}) \}$ be a dominant 
sequence of rectangular shape partitions. Denote by $\lambda^{'}$ the conjugate of
 $\ld$, and by $R^{'}$ a dominant rearrangement of a sequence of rectangular shape 
partitions $\{ (\eta_{a}^{\mu_{a}}) \}$. Then
$$ K_{\ld,R}(q) = q^{n(R)} K_{\ld^{'},R^{'}}(q^{-1}),$$
where
$$n(R)= \sum_{a < b} \min(\mu_{a},\mu_{b}) \min(\eta_{a},\eta_{b}).$$ 
\end{theorem}

A technical proof is based on checking of the statement that the map
$$ \iota:  m_{ij} \longrightarrow {\hat{m}}_{ij} =-m_{ji}+\theta(\ld_j-i)+\sum_{a\ge 1}\theta(\mu_a-j) \theta(\eta_a-i)$$
establishes bijection between the sets of admissible configurations of types 
$(\ld, R)$ and $(\ld^{'},R^{'})$, and $\iota(c(m_{ij})) =c(({\hat{m}}_{ij})).$ 

\subsection{Example} {\rm

 Let $n=6$, consider for example, a  standard Young tableau 
$$ T = \begin{array}{cccccc} 1 & 2 & 3 & 6 & 8 & 9\\
4 & 5 & 7 & 1 0 & 1 1 & 12 \end{array},~~~c(T)=48.$$
The corresponding rigged configuration $(\nu,J)$ is 
$$\nu =(321),~~J=(J_3 =0,~J_2=2, J_1 = 6),~~(m_{ij})(\nu) = 
 \left( \begin{array}{ccc} 9 & -2 & -1 \\
3 & 2 & 1 \end{array} \right),~~~c(\nu)=44. $$
 Recall that $c(T)$ and $c(\nu)$ denote the charge of tableau $T$ and configuration $\nu$ correspondingly.

$\bullet$~~One can see that $c(T)=c(\nu)+J_3+J_2+J_1$, as it should be in general. 

$\bullet$~~Now, the descent set and  descent number of tableau $T$ are $Des(T)=\{3,6,9 \}$,~$des(T)=3$. One can see that $des(T)= 3= \nu_1^{'}$, as it should be in general
\footnote{~~ In fact the shape of the first configuration $\nu^{(1)}$ of type
 $(\ld,\mu)$ can be read off from the set of ~''secondary''~~descent sets 
$ \{ Des^{(1)} (T)=Des(T), Des^{(2)}(T), \ldots,...)$, cf  \cite{Ki10}.}. 

$\bullet$~~~One can check that our tableau $T$ is invariant under the action 
of the Sch\"{u}tzenberger involution
\footnote{~~$http://en.wikipedia.org/wiki/Jeu\_de\_taquin$ } 
on the set of standard Young tableaux of a shape  $\ld$. It is clearly seen 
from  the  set of riggings $J$
\footnote{~~In our example ~$J=(0,1,3)$.}
that the rigged configuration $(\nu, J)$ corresponding to tableau $T$,  is 
invariant under the Flip involution
\footnote{~~Recall that a rigging of an admissible configuration $\nu$ is a 
collection of integers 
$$J=(\{J_{s,r}^{(k)}\}~~,~1\le s \le m_{r}(\nu^{(k)})$$ 
 such  that for a given $k,r$  one has
$$ 0 \le J_{1,r}^{(k)} \le J_{2,r}^{(k)} \le \ldots \le 
J_{{m_{r}(\nu^{(k)}),r}}^{(k)} \le P_{r}^{(k)}(\nu).$$
The Flip involution $\kappa$ is defined as follows:
$$\kappa(\nu, \{J_{s,r}^{(k)}\}) =(\nu,\{ J_{m_{r}(\nu^{(k)}) - s +1,r}^{(k)}\}).$$
}
 on the set of rigged configurations of type $(\ld, 1^{|\ld|})$, as it should 
be in general, see \cite{KSS} for a complete proof of the statement that the 
action of the Sch\"{u}tzenberger transformation  on a Littlewood--
Richardson tableau  $T \in LR(\ld, R)$,  under the Rigged Configuration 
Bijection transforms tableau $T$ to a Littlewood--Richardson tableau 
corresponding to  the rigged configuration $\nu \kappa(J))$, where $(\nu J)$ 
is the rigged configuration corresponding to tableau T we are started with.

}

% ${\s}_w$
\end{document}